\documentclass[english,nofootinbib, notitlepage, showkeys, final]{revtex4-1}
\usepackage[T1]{fontenc}
\usepackage[latin9]{inputenc}
\usepackage{float}
\usepackage{amsmath}
\usepackage{amssymb}

\makeatletter

\newcommand{\noun}[1]{\textsc{#1}}
\providecommand{\tabularnewline}{\\}

\usepackage{tikz}
\usepackage{pst-tree}
\usepackage{tikz} 
\usepackage{tikz-qtree} 
\usetikzlibrary{calc,positioning,shadows,arrows,trees,shapes,snakes,fit,intersections}
\allowdisplaybreaks
\usepackage{bm}
\usepackage[final]{microtype}
\usepackage[unicode=true, bookmarks=true,bookmarksnumbered=false,
bookmarksopen=false, breaklinks=true,pdfborder={0 0 0},
pdfborderstyle={},backref=false,colorlinks=true,citecolor=blue]{hyperref}

\AtBeginDocument{
  
}

\makeatother

\usepackage{babel}
\begin{document}

\title{The algebraic Bethe Ansatz and combinatorial trees}

\author{R. S. Vieira }
\email{rsvieira@df.ufscar.br}

\address{Universidade Estadual Paulista, Faculdade de Ciências e Tecnologia,
Departamento de Matemática e Computação, Caixa-Postal 467, CEP. 19060-900,
Presidente Prudente, SP, Brasil}

\author{A. Lima-Santos}
\email{dals@df.ufscar.br}

\address{Universidade Federal de São Carlos, Departamento de Física, Caixa-Postal
676, CEP. 13565-905, São Carlos, SP, Brasil.}
\begin{abstract}
We present in this paper a comprehensive introduction to the algebraic
Bethe Ansatz, taking as examples the six-vertex model with periodic
and non-periodic boundary conditions. We propose a diagrammatic representation
of the commutation relations used in the algebraic Bethe Ansatz, so
that the action of the transfer matrix in the $n$th excited state
gives place to labeled combinatorial trees. The analysis of these
combinatorial trees provides in a straightforward way the eigenvalues
and eigenstates of the transfer matrix, as well as the respective
Bethe Ansatz equations. Several identities between the $R$-matrix
elements can also be derived from the symmetry of these diagrams regarding
the permutation of their labels. This combinatorial approach gives
some insights about how the algebraic Bethe Ansatz works, which can
be valuable for non-experts readers.
\end{abstract}

\maketitle

\section{Introduction}

The algebraic Bethe Ansatz is a powerful technique for solving analytically
the eigenvalue problem in many-body quantum field theory and statistical
mechanics. This method was created by the Leningrad group in the context
of quantum field theory, as a quantum generalization of the inverse
scattering method \cite{Faddeev1995}. Soon after, the same mathematical
structure led to an algebraic formulation of the (coordinate) Bethe
Ansatz \cite{Faddeev1995B}, this time in the field of statistical
mechanics. Several models of quantum field theory and statistical
mechanics were successfully solved by the algebraic Bethe Ansatz \cite{Korepin1997}:
the non-linear Schrödinger equation \cite{Sklyanin1979A} and the
Sine-Gordon model \cite{Sklyanin1982} are typical examples in quantum
field theory; in statistical mechanics we can cite the one-dimensional
\noun{xxx}, \noun{xxz} and \noun{xyz} Heisenberg spin chains, the
six-vertex model and the eight-vertex model \cite{TakhtadzhanFaddeev1979}. 

Although the algebraic Bethe Ansatz is by now a very well understood
method, it is nevertheless very technical and usually a difficult
matter for students or non-experts researches. In this paper we present
a comprehensive introduction to the algebraic Bethe Ansatz, taking
as example the six-vertex model with both periodic and non-periodic
boundary conditions. We also propose \textendash{} which is the main
novelty of the paper \textendash{} a diagrammatic representation for
the commutation relations used in the algebraic Bethe Ansatz that
provides the eigenvalues and eigenstates of the transfer matrix, and
also the Bethe Ansatz equations, in a straightforward way. In fact,
in the usual algebraic Bethe Ansatz, one needs to use the aforementioned
commutation relations repeatedly in order to compute the action of
the transfer matrix on the excited states, which is generally very
cumbersome. In our diagrammatic approach, however, this action is
represented by simple combinatorial trees and the results follow through
a combinatorial analysis only.

This paper is organized as follows: in section \noun{\ref{Section-Periodic}}
we consider the six-vertex model with periodic boundary conditions.
We first introduce some statistical concepts as the monodromy and
transfer matrices in section \ref{Section-SM-P} and, then, we explain
step-by-step in section \ref{Section-ABA-P} how the algebraic Bethe
Ansatz works. We present our combinatorial approach in section \ref{Section-Combinatorial-P},
where we show how the eigenvalues and eigenstates of the transfer
matrix \textendash{} and also the Bethe Ansatz equations \textendash{}
can be easily derived from the analysis of simple diagrams. In section
\ref{Section-Boundary}, we study the six-vertex model with non-periodic
boundary conditions. We discuss the role of the boundaries on the
monodromy and transfer matrices in section \ref{Section-SM-B} and,
in section \ref{Section-ABA-B}, we also give a step-by-step presentation
of the boundary algebraic Bethe Ansatz. The combinatorial approach
for the non-periodic case is presented in \ref{Section-Combinatorial-B}
and we close the paper with a brief conclusion in section \ref{Section-Conclusion}. 

\section{The six-vertex model with periodic boundary conditions\label{Section-Periodic}}

\subsection{The monodromy and transfer matrices of the six-vertex model with
periodic boundary conditions \label{Section-SM-P}}

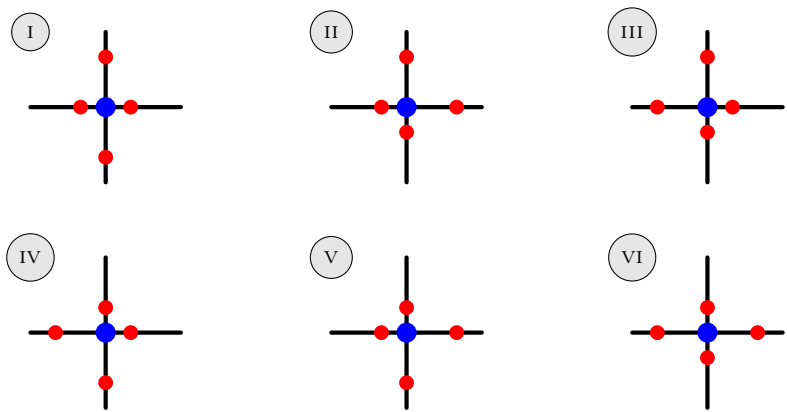
\begin{figure}[H]
\begin{center}
\begin{tikzpicture}[ultra thick, line cap=round, roundnode/.style={circle, draw=black, fill=black!10, thin, minimum size=5mm}]

\draw (-5,0) -- (-3,0);\draw (-4,-1) -- (-4,1);
\draw (-1,0) -- (+1,0);\draw (+0,-1) -- (+0,1);
\draw (+3,0) -- (+5,0);\draw (+4,-1) -- (+4,1);

\draw (-5,-3) -- (-3,-3);\draw (-4,-4) -- (-4,-2);
\draw (-1,-3) -- (+1,-3);\draw (+0,-4) -- (+0,-2);
\draw (+3,-3) -- (+5,-3);\draw (+4,-4) -- (+4,-2);

\filldraw[blue] (-4 ,0   ) circle (3pt);
\filldraw[red] (-4-1/3 ,0   ) circle (2pt);
\filldraw[red] (-4+1/3 ,0   ) circle (2pt);
\filldraw[red] (-4     ,+2/3) circle (2pt);
\filldraw[red] (-4     ,-2/3) circle (2pt);

\filldraw[blue] (0 ,0   ) circle (3pt);
\filldraw[red] (-0-1/3 ,0   ) circle (2pt);
\filldraw[red] (-0+2/3 ,0   ) circle (2pt);
\filldraw[red] (-0     ,+2/3) circle (2pt);
\filldraw[red] (-0     ,-1/3) circle (2pt);

\filldraw[blue] (4 ,0   ) circle (3pt);
\filldraw[red] (+4-2/3 ,0   ) circle (2pt);
\filldraw[red] (+4+1/3 ,0   ) circle (2pt);
\filldraw[red] (+4     ,+2/3) circle (2pt);
\filldraw[red] (+4     ,-1/3) circle (2pt);

\filldraw[blue] (-4 ,-3   ) circle (3pt);
\filldraw[red] (-4-2/3 ,-3   ) circle  (2pt);
\filldraw[red] (-4+1/3 ,-3   ) circle  (2pt);
\filldraw[red] (-4     ,-3 +1/3) circle (2pt);
\filldraw[red] (-4     ,-3 -2/3) circle (2pt);

\filldraw[blue] (0 ,-3   ) circle (3pt);
\filldraw[red] (-0-1/3 ,-3   ) circle  (2pt);
\filldraw[red] (-0+2/3 ,-3   ) circle  (2pt);
\filldraw[red] (-0     ,-3+1/3) circle (2pt);
\filldraw[red] (-0     ,-3-2/3) circle (2pt);

\filldraw[blue] (4 ,-3   ) circle (3pt);
\filldraw[red] (+4-2/3 ,-3   ) circle  (2pt);
\filldraw[red] (+4+2/3 ,-3   ) circle  (2pt);
\filldraw[red] (+4     ,-3+1/3) circle (2pt);
\filldraw[red] (+4     ,-3-1/3) circle (2pt);

\coordinate (A) at (-5,1);
\coordinate (B) at (-1,1);
\coordinate (C) at (+3,1);
\coordinate (D) at (-5,-2);
\coordinate (E) at (-1,-2);
\coordinate (F) at (+3,-2);

\node[roundnode] at (A) {\textsc{i}};
\node[roundnode] at (B) {\textsc{ii}};
\node[roundnode] at (C) {\textsc{iii}};
\node[roundnode] at (D) {\textsc{iv}};
\node[roundnode] at (E) {\textsc{v}};
\node[roundnode] at (F) {\textsc{vi}};

\end{tikzpicture} 
\end{center}

\caption{The six-vertex configurations of the square ice model, according to
the Pauling ice rule. The oxygen atoms are situated in the center
of the vertices, while the hydrogen atoms are on the edges. In each
edge there is exactly one hydrogen atom, which can be near or far
from the oxygen atom. The Pauling ice rule states that the only physical
configurations are those in which there is two, and only two, hydrogen
atoms close to any oxygen atom. This rule agrees with the usual molecular
formula for the water: \noun{h$_{2}$o.}}

\label{Fig6V}
\end{figure}

The first model to be solved through the algebraic Bethe Ansatz were
the so-called \emph{six-vertex model}. This was a statistical model
introduced by Linus Pauling \cite{Pauling1935} in an attempt to explain
the statistical properties of the water ice \textendash{} more specifically,
to get account for the residual entropy of the ice. Remember that
at low temperatures, the molecules of water arrange into an almost
crystalline lattice. Pauling considered a two-dimensional approximation
for the water ice lattice in which each oxygen atom is disposed in
a \emph{vertex} of a square lattice and that there is a hydrogen atom
on each edge of this vertex. Each hydrogen atom is supposed to be
either near or far from the oxygen atom, so that we have in total
sixteen vertex configurations \textendash{} which lead us to a \emph{sixteen
vertex model} \textendash{} and for each vertex configuration, a respective
energy and a Boltzmann weight is associated. However, the fact that
each water molecule has indeed two, and only two, bound hydrogens
implies that some configurations should be despised. This restriction
was already noticed by Pauling in \cite{Pauling1935} and, since then,
it is known as \emph{Pauling ice rule}. Taking into account the Pauling
ice rule, the number of possible configurations of a given vertex
representing a water molecule reduces to only six, and we get a \emph{six-vertex
model}. These six allowed configurations are illustrated in Fig. \ref{Fig6V}.

Now, let us see how we can construct the partition function for the
six-vertex model. To this end, consider a square lattice with $L$
columns and $N$ lines. We can impose periodic boundary conditions
in both the horizontal and vertical directions, which means that any
vertex at the position $(i+N,j+L)$ is to be identified with the vertex
at the position $(i,j)$ of the lattice. To each vertex at the position
$(i,j)$ of the lattice we associate a horizontal local Hilbert space
$H_{i}$ $(1\leqslant i\leqslant N)$ and a vertical local Hilbert
space $V_{j}$ $(1\leqslant j\leqslant L)$, both isomorphic to $\mathbb{C}^{2}$,
so that the local Hilbert space of the vertex can be written as $H_{i,j}=H_{i}\otimes V_{j}$,
which is isomorphic to $\mathbb{C}^{2}\otimes\mathbb{C}^{2}$. Therefore,
the Hilbert space associated with the whole lattice can be written
as $\mathcal{H}=H\otimes V$, where $H=H_{1}\otimes\cdots\otimes H_{L}$
stands for the horizontal spaces and $V=V_{1}\otimes\cdots\otimes V_{N}$
for the vertical ones. Notice that $\mathcal{H}$ is isomorphic to
$\mathbb{C}^{2N}\otimes\mathbb{C}^{2L}$. 

For a given vertex at the position $(i,j)$ of the lattice, we introduce
an $R$-matrix whose elements are related to the Boltzmann weights
associated with the possible configurations of the vertices. The $R$-matrix
corresponding to a vertex at the point $(i,j)$ of the lattice has
values in $\mathrm{End}\left(H_{i}\otimes V_{j}\right)$, so that
it is a four-by-four matrix. Since there are only six possible vertex
configurations, the $R$-matrix must have only six non-null entries.
Besides, if we take into account some physical symmetries (e.g. time-inversion,
parity etc.), then some vertex configurations should be equivalent,
which means that some elements of the $R$-matrix should be the same.
It turns out then that the $R$-matrix for the six-vertex model can
be written as, 
\begin{equation}
R(u)=\begin{pmatrix}r_{1}(u) & 0 & 0 & 0\\
0 & r_{2}(u) & r_{3}(u) & 0\\
0 & r_{3}(u) & r_{2}(u) & 0\\
0 & 0 & 0 & r_{1}(u)
\end{pmatrix},
\end{equation}
where the amplitudes $r_{1}(u)$, $r_{2}(u)$ and $r_{3}(u)$ are
related to the Boltzmann weights of the vertex configurations (the
exact expressions for them will be presented below) and $u$ is the
so-called \emph{spectral parameter.}

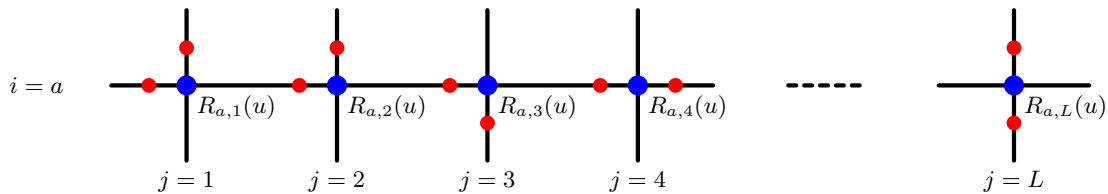
\begin{figure}
\begin{center}
\begin{tikzpicture}[ultra thick,line cap=round]

\draw (0,0) -- (8,0);
\draw[dashed] (9,0) -- (10,0);
\draw (11,0) -- (13,0);

\draw (1,-1) -- (1,1);
\draw (3,-1) -- (3,1);
\draw (5,-1) -- (5,1);
\draw (7,-1) -- (7,1);
\draw (12,-1) -- (12,1);

\filldraw[blue] (1 ,0) circle  (3pt);
\filldraw[blue] (3 ,0) circle  (3pt);
\filldraw[blue] (5 ,0) circle  (3pt);
\filldraw[blue] (7 ,0) circle  (3pt);
\filldraw[blue] (12 ,0) circle (3pt);

\filldraw[red] (0.5 ,0) circle  (2pt);
\filldraw[red] (2.5 ,0) circle  (2pt);
\filldraw[red] (4.5 ,0) circle  (2pt);
\filldraw[red] (6.5 ,0) circle  (2pt);
\filldraw[red] (7.5 ,0) circle  (2pt);

\filldraw[red] (1 , 0.5) circle  (2pt);
\filldraw[red] (3 , 0.5) circle  (2pt);
\filldraw[red] (12, 0.5) circle  (2pt);

\filldraw[red] (5 , -0.5) circle  (2pt);
\filldraw[red] (12, -0.5) circle  (2pt);

\node[left]  at (-0.5,0) {$i=a$};
\node[below] at (1,-1) {$j=1$};
\node[below] at (3,-1) {$j=2$};
\node[below] at (5,-1) {$j=3$};
\node[below] at (7,-1) {$j=4$};
\node[below] at (12,-1) {$j=L$};

\node[below right] at (1,0) {$R_{a,1}(u)$};
\node[below right] at (3,0) {$R_{a,2}(u)$};
\node[below right] at (5,0) {$R_{a,3}(u)$};
\node[below right] at (7,0) {$R_{a,4}(u)$};
\node[below right] at (12,0) {$R_{a,L}(u)$};

\end{tikzpicture} 
\end{center}

\caption{A row of the lattice. The blue circles at the center of the vertices
represent the oxygen atoms. The red circles at the edges represent
the surrounding hydrogen atoms. Notice that each oxygen atom has two
orbiting hydrogens, so that the Pauling ice rule is satisfied. The
energy $E$ of each water molecule, and hence the Boltzmann weights
$R$, depends only on the corresponding vertex configurations, that
is, regarding the distance \textendash{} near or far \textendash{}
of the two hydrogen atoms with respect to the oxygen atom. The monodromy
matrix is given by the product of all $R$-matrix in a row of the
lattice. The transfer matrix is found after we take the trace of the
monodromy matrix in the horizontal vector space. }

\label{Fig-Row}
\end{figure}

Now, let us take a given row of the lattice, say the row $i=a$. Taking
the product of all $R$-matrices in this row, we get the so-called
\emph{monodromy matri}x\footnote{The monodromy matrix $M_{a}(u)$ acts in $\mathrm{End}\left(H_{a}\otimes V_{1}\otimes\cdots\otimes V_{L}\right)$.
By this reason, the $R_{aq}(u)$ matrices that appears in the definition
of $M_{a}(u)$ should be regarded as matrices with values $\mathrm{End}\left(H_{a}\otimes V_{1}\otimes\cdots\otimes V_{L}\right)$
that act non-trivially only in $\mathrm{End}\left(H_{a}\otimes V_{q}\right)$
and as the identity matrix in the other remaining vector spaces. The
space $H_{a}$ is usually called \emph{auxiliary space}, while the
spaces $V_{q}$ $(1\leqslant q\leqslant L)$ are called \emph{quantum
spaces}.}, 
\begin{equation}
M_{a}(u)=R_{a1}(u)R_{a2}(u)\cdots R_{aL}(u),\label{M}
\end{equation}
whose elements are related to the Boltzmann weights associated with
all possible configurations of the row \textendash{} see Fig. \ref{Fig-Row}.
The trace of the monodromy matrix on the space $H_{a}$ provides a
sum over all possible row configurations, that is, it gives a reduced
partition function for this row, which is called \emph{transfer matrix}:
\begin{equation}
T_{a}(u)=\text{tr}_{a}\left[M_{a}(u)\right]=\text{tr}_{a}\left[R_{a1}(u)R_{a2}(u)\cdots R_{aL}(u)\right].
\end{equation}
From this, we can easily find the total partition function of the
system. In fact, this follows after we multiply all the monodromy
matrices of each lattice row and we then take the trace:
\begin{equation}
Z=\text{tr}\left[M_{1}(u)M_{2}(u)\cdots M_{N}(u)\right].
\end{equation}

Notice, moreover, that when the rows of the lattice are all equivalent,
the eigenvalues of the transfer matrices associated with any row of
the lattice will be the same. The, if $\tau_{1},\ldots,\tau_{L}$
denote the eigenvalues of the transfer matrix, the corresponding eigenvalues
of the partition function will be,
\begin{equation}
Z=\tau_{1}^{N}+\cdots+\tau_{L}^{N}.
\end{equation}
This means that we can look for the diagonalization of the transfer
matrix only. This is exactly what the algebraic Bethe Ansatz concerns
with, as we shall see in the following sections.

\subsection{The algebraic Bethe Ansatz\label{Section-ABA-P}}

The algebraic Bethe Ansatz provides an analytical solution for the
diagonalization of the transfer matrix for integrable models as the
six-vertex model with periodic boundary conditions. It is implemented,
however, in several steps which we may enumerate as follows:
\begin{enumerate}
\item \emph{A solution of the Yang-Baxter equation providing the $R$-matrix}.
The starting point is a given $R$-matrix, solution of the Yang-Baxter
equation, that describes the model. The Yang-Baxter equation ensures
the integrability of the model, which means that it can be solved
in an exact way.
\item \emph{The Lax representation of the monodromy and transfer matrices.
}This is a representation in which the monodromy explicitly exhibit
the annihilator and creator operators used in the construction of
the excited states. Besides, the transfer matrix become given by sum
of the diagonal monodromy elements in this representation;
\item \emph{The reference state.} It is just a simple enough eigenstate
of the transfer matrix in which the corresponding eigenvalue can be
evaluated directly; 
\item \emph{The construction of the excited states}. They are built through
the action of the aforementioned creator operators on the reference
state;l
\item \emph{The derivation of the commutation relations}. The commutation
relations are necessary to one compute the action of the transfer
matrix on the excited states;
\item \emph{The computation of the eigenvalues} of the transfer matrix.
This is the most difficult step. To compute the eigenvalues we need
to use the commutation relations mentioned above repeatedly;
\item \emph{The solution of the Bethe Ansatz equations}. These are a system
of coupled non-linear equations that appears as consistency conditions
of the method. The Bethe Ansatz equations need to be solved in order
to one obtain a completely analytical solution for the spectral problem. 
\end{enumerate}
Next we shall explain in details how the steps above are developed.

\subsubsection{The $R$-matrix, solution of the Yang-Baxter equation}

The start point of the periodic algebraic Bethe Ansatz is the existence
a given $R$-matrix, a solution of the Yang-Baxter equation \cite{Yang1967,Baxter1972},
\begin{equation}
R_{12}(u-v)R_{13}(u)R_{23}(v)=R_{23}(v)R_{13}(u)R_{12}(u-v).\label{YBE}
\end{equation}

This is a matrix equation with values in $\mathrm{End}\left(V_{a}\otimes V_{b}\otimes V_{c}\right)$,
where $V_{a}$, $V_{b}$ and $V_{c}$ are complex vector spaces isomorphic
to $\mathbb{C}^{2}$ (in the case of the six-vertex model). The operators
$R_{ab}$ act as an $R$-matrix in $\mathrm{End}\left(V_{a}\otimes V_{b}\right)$
and as the identity in the other vector space $V_{c}$. 

\begin{figure}[H]
\begin{center}
\begin{tikzpicture}[ultra thick,line cap=round]
\draw (-1.5,-1) -- (-1, -1);
\draw (-1.5, 1) -- (-1,  1);
\draw (-1,-1) -- (1, 1);
\draw (-1, 1) -- (1,-1);
\draw ( 1, 1) -- (3, 1);
\draw ( 1,-1) -- (3,-1);
\draw ( 2,-2) -- (2, 2);
\node at (4,0){$\bm{=}$};
\draw ( 5, 1) -- (7, 1);
\draw ( 5,-1) -- (7,-1);
\draw ( 6,-2) -- (6, 2);
\draw ( 7,-1) -- (9, 1);
\draw ( 7, 1) -- (9,-1);
\draw (9 ,-1) -- (9.5, -1);
\draw (9 , 1) -- (9.5,  1);

\filldraw[red]  (0, 0) circle  (3pt);
\filldraw[blue] (2,-1) circle  (3pt);
\filldraw[blue] (2, 1) circle  (3pt);
\filldraw[blue] (6, 1) circle  (3pt);
\filldraw[blue] (6,-1) circle  (3pt);
\filldraw[red]  (8, 0) circle  (3pt);

\node[left]       at (-0.3, 0){$S_{12}(u)$};
\node[below left] at (   2,-1){$R_{13}(u)$};
\node[above left] at (   2, 1){$R_{23}(u)$};

\node[right]        at (8.3, 0){$S_{12}(u)$};
\node[below right] at (  6,-1){$R_{13}(u)$};
\node[above right] at (  6, 1){$R_{23}(u)$};

\node[below] at (  2,-2){$3$};
\node[below] at (  6,-2){$3$};
\node[right] at (  3,-1){$1$};
\node[left]  at (  5,-1){$1$};
\node[right] at (  3, 1){$2$};
\node[left]  at (  5, 1){$2$};

\end{tikzpicture} 
\end{center}

\caption{A graphical representation of the Yang-Baxter Equation (\ref{YBE-SR}).
The $S$-matrix ``twist'' two adjacent rows of the lattice. The Yang-Baxter
equation states that if we twist the rows $1$ and $2$ and then we
take the product of the $R$ matrix at the vertex $(1,3)$ and $(2,3)$,
respectively, then we shall get the same result if we multiply these
$R$ matrices first and then we twist the rows. }

\label{Fig-YB}
\end{figure}
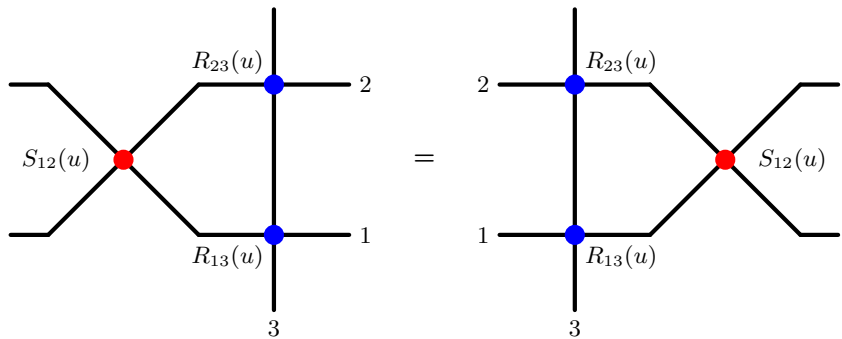

Introducing the permutator matrix, 
\begin{equation}
P=\begin{pmatrix}1 & 0 & 0 & 0\\
0 & 0 & 1 & 0\\
0 & 1 & 0 & 0\\
0 & 0 & 0 & 1
\end{pmatrix},
\end{equation}
 whose action on any element $A\otimes B\in\mathrm{End}\left(V_{a}\otimes V_{b}\right)$
is $P\left(A\otimes B\right)P=B\otimes A$, we can define a ``twisted
$R$-matrix'', also called $S$-matrix, as follows: 
\begin{equation}
S(u)=PR(u)=\begin{pmatrix}r_{1}(u) & 0 & 0 & 0\\
0 & r_{3}(u) & r_{2}(u) & 0\\
0 & r_{2}(u) & r_{3}(u) & 0\\
0 & 0 & 0 & r_{1}(u)
\end{pmatrix}.\label{S}
\end{equation}
 This matrix allows us to rewrite the Yang-Baxter equation in a more
symmetric form, namely, as 
\begin{equation}
S_{12}(u-v)S_{23}(u)S_{12}(v)=S_{23}(v)S_{12}(u)S_{23}(u-v).\label{YBE-S}
\end{equation}
 There is also a mixed representation that is useful:
\begin{equation}
S_{12}(u-v)R_{13}(u)R_{23}(v)=R_{13}(v)R_{23}(u)S_{12}(u-v),\label{YBE-SR}
\end{equation}
which follows from the properties of the permutator matrix (e.g.,
the identity $P_{ab}R_{bc}R_{ac}P_{ab}=R_{ac}R_{bc}$). A graphical
representation of the Yang-Baxter equation is given in Fig. \ref{Fig-YB}.

Now, let us see how we can find the solutions of the Yang-Baxter equation
for the six-vertex model $R$-matrix. From what we have seen on the
previous section, the $R$-matrix we are looking for should have the
following form\footnote{We had assumed the most symmetrical case, which corresponds to the
full symmetric six-vertex model. There are solutions of the Yang-Baxter
equations for the non-symmetric cases as well. We indicate the reference
\cite{Vieira2017b} for the derivation of these most general solutions.}: 
\begin{equation}
R(u)=\begin{pmatrix}r_{1}(u) & 0 & 0 & 0\\
0 & r_{2}(u) & r_{3}(u) & 0\\
0 & r_{3}(u) & r_{2}(u) & 0\\
0 & 0 & 0 & r_{1}(u)
\end{pmatrix}.\label{R0}
\end{equation}
Regardless what form of the Yang-Baxter equation we use, they represent
a system of coupled functional equations for the $R$-matrix elements.
Several equations, however are dependent on each other or are null,
so that we have actually only three independent equations, namely,
\begin{align}
r_{1}(u-v)r_{1}(u)r_{3}(v) & =r_{2}(u-v)r_{2}(u)r_{3}(v)+r_{3}(u-v)r_{3}(u)r_{2}(v),\\
r_{3}(u-v)r_{1}(u)r_{2}(v) & =r_{3}(u-v)r_{2}(u)r_{1}(v)+r_{2}(u-v)r_{3}(u)r_{3}(v),\\
r_{1}(u-v)r_{3}(u)r_{2}(v) & =r_{2}(u-v)r_{3}(u)r_{1}(v)+r_{3}(u-v)r_{2}(u)r_{3}(v).
\end{align}
 This system of functional equations can be solved through a kind
of separation of variables method \cite{Baxter2016}. Eliminating
the variables with the dependence on $u-v$, we shall obtain the separated
equation, 
\begin{equation}
\frac{r_{1}^{2}(u)+r_{2}^{2}(u)+r_{3}^{2}(u)}{r_{1}(u)r_{2}(u)}=\frac{r_{1}^{2}(v)+r_{2}^{2}(v)+r_{3}^{2}(v)}{r_{1}(v)r_{2}(v)}=\varDelta,\label{Delta}
\end{equation}
where $\varDelta$ is a constant independent of $u$ and $v$. The
solutions of (\ref{YBE}) follow after we find a parameterization
of $r_{1}(u)$, $r_{2}(u)$ and $r_{3}(u)$ satisfying (\ref{Delta}).
It follows that there are several possible solutions of this functional
equation, each solution referring to a specific six-vertex model.
The main important ones are the solution 
\begin{equation}
r_{1}(u)=\sinh\left(u+\xi\right),\qquad r_{2}(u)=\sinh u,\qquad r_{3}(u)=\sinh\xi,\label{Relements}
\end{equation}
 which is related to the \textsc{xxz} Heisenberg chain, and the solution
\begin{equation}
r_{1}(u)=u+\xi,\qquad r_{2}(u)=u,\qquad r_{3}(u)=\xi,
\end{equation}
which is related to the \textsc{xxx} Heisenberg chain (in these expressions,
$\xi$ is an arbitrary parameter). 

\subsubsection{The Lax representation of the monodromy and transfer matrices}

Once we have an $R$-matrix, solution of the Yang-Baxter equation,
we pass to construct the monodromy and transfer matrices. We conveniently
use the a representation \textendash{} the \emph{Lax representation}
\textendash{} in which the $R$-matrix (\ref{R0}) is written as a
two-by-two operator-valued matrix:
\begin{equation}
R(u)=\begin{pmatrix}L_{1}^{1}(u) & L_{1}^{2}(u)\\
L_{2}^{1}(u) & L_{2}^{2}(u)
\end{pmatrix},\label{R}
\end{equation}
where, 
\begin{align}
L_{1}^{1}(u) & =\begin{pmatrix}r_{1}(u) & 0\\
0 & r_{2}(u)
\end{pmatrix}, & L_{1}^{2}(u) & =\begin{pmatrix}0 & 0\\
r_{3}(u) & 0
\end{pmatrix}, & L_{2}^{1}(u) & =\begin{pmatrix}0 & r_{3}(u)\\
0 & 0
\end{pmatrix}, & L_{2}^{2}(u) & =\begin{pmatrix}r_{2}(u) & 0\\
0 & r_{1}(u)
\end{pmatrix}.\label{Lij}
\end{align}

As we have seen in the previous section, the \emph{monodromy matrix}
is given by the product of the $R$-matrices running through all the
sites of a given row the lattice. That is, for a given $i=a$ of the
lattice, we have, according to (\ref{M}), 
\begin{equation}
M_{a}(u)=R_{a1}(u)\cdots R_{aL}(u).
\end{equation}
Here, $R_{aq}$ $(1\leqslant q\leqslant L)$ means an $R$-matrix
that act non-trivially only in $\mathrm{End}\left(H_{a}\otimes V_{q}\right)$.
In the Lax representation, however, the monodromy matrix becomes a
two-by-two matrix with values in $\mathrm{End}\left(H_{a}\right)$
\textendash{} its elements are operators acting in $\mathrm{End}\left(V_{1}\otimes\cdots\otimes V_{L}\right)$:
\begin{equation}
M_{a}(u)=\begin{pmatrix}M_{1}^{1}(u) & M_{1}^{2}(u)\\
M_{2}^{1}(u) & M_{2}^{2}(u)
\end{pmatrix}\equiv\begin{pmatrix}A(u) & B(u)\\
C(u) & D(u)
\end{pmatrix}.\label{MonoP}
\end{equation}
 In terms of the operators (\ref{Lij}), the monodromy elements $M_{i}^{j}(u)$
can be found explicitly by the formula: 
\begin{equation}
M_{i}^{j}(u)=\sum_{k_{1},\ldots,k_{L-1}=1}^{2}L_{i}^{k_{1}}(u)\otimes L_{k_{1}}^{k_{2}}(u)\otimes\cdots\otimes L_{k_{L-2}}^{k_{L-1}}(u)\otimes L_{k_{L-1}}^{j}(u).\label{Mij}
\end{equation}

Finally, the \emph{transfer matrix} is defined as the trace (in the
$H_{a}$ space) of the monodromy matrix, that is,
\begin{equation}
T(u)=\mathrm{tr}_{a}\left[M_{a}(u)\right]=A(u)+D(u).\label{tP}
\end{equation}

\subsubsection{The reference state}

The next step in the execution of the algebraic Bethe Ansatz consists
of finding an appropriate reference state. This corresponds to a (simple
enough) eigenstate of the transfer matrix so that its eigenvalue can
be directly computed. 

Fortunately, we can verify that the reference state for the six-vertex
model is the most simple possible state, namely,
\begin{equation}
\Psi_{0}=\begin{pmatrix}1\\
0
\end{pmatrix}_{1}\otimes\cdots\otimes\begin{pmatrix}1\\
0
\end{pmatrix}_{L}.
\end{equation}
In fact, it follows from (\ref{Lij}) and (\ref{Mij}) that the action
of the monodromy elements on the reference state is given by
\begin{equation}
A(u)\Psi_{0}=\alpha(u)\Psi_{0},\qquad B(u)\Psi_{0}\neq z\Psi_{0},\qquad C(u)\Psi_{0}=0\Psi_{0},\qquad D(u)\Psi_{0}=\delta(u)\Psi_{0},
\end{equation}
 where $z$ can be any complex number and
\begin{equation}
\alpha(u)=r_{1}^{L}(u),\qquad\delta(u)=r_{2}^{L}(u).
\end{equation}
From this is not difficult to show that the action of the transfer
matrix on $\Psi_{0}$ is given by, 
\begin{equation}
T(u)\Psi_{0}=\tau_{0}(u)\Psi_{0},\qquad\mathrm{where},\qquad\tau_{0}(u)=\alpha(u)+\delta(u),
\end{equation}
so that $\Psi_{0}$ is indeed a eigenstate of the transfer matrix
whose eigenvalue is $\tau_{0}(u)$.

Notice moreover that the operator $C(u)$ annihilates $\Psi_{0}$,
while the operator $B(u)$ gives something not proportional to $\Psi_{0}$
\textendash{} that is it creates another state. By this reason we
say that the $C$ operators are \emph{annihilator operators}, while
the $B$ operators are \emph{creator operators}. We shall call the
$A$ and $D$ operators as \emph{diagonal operators}. 

\subsubsection{The construction of the excited states}

Once the reference state is determined, we can construct the \emph{excited
states} of the transfer matrix by acting with the creator operator
$B$ repeatedly on the reference state $\Psi_{0}$. In this way, we
define $n$th excited state of the transfer matrix as
\begin{equation}
\Psi_{n}\left(u_{1},\ldots,u_{n}\right)=B(u_{1})\cdots B(u_{n})\Psi_{0}.\label{PsiNP}
\end{equation}
Notice that each creator operator in (\ref{PsiNP}) depends on a different
parameter $u_{k}$ $(1\leqslant k\leqslant n)$; these parameters
are called \emph{rapidities} and they are until now undetermined.
We hope to fix the values of the rapidities so that $\Psi_{n}$ becomes
indeed an eigenstate of the transfer matrix. We shall see in a while
that the rapidities can be fixed (at least implicitly) by a system
of non-linear equations called the \emph{Bethe Ansatz equations}. 

\subsubsection{The commutation relations}

Our go now is to evaluate the action of the transfer matrix on these
excited states. From (\ref{tP}) and (\ref{PsiNP}), we get that,
\begin{equation}
T(u)\Psi_{n}\left(u_{1},\ldots,u_{n}\right)=A(u)B(u_{1})\cdots B(u_{n})\Psi_{0}+D(u)B(u_{1})\cdots B(u_{n})\Psi_{0}.\label{TPsi6VP}
\end{equation}

\begin{figure}
\begin{center}
\begin{tikzpicture}[ultra thick,line cap=round]

\draw (-0.5,0) -- (0.5 ,2);
\draw (-0.5,2) -- (0.5 ,0);
\draw (-1  ,2) -- (-0.5,2);
\draw (-1  ,0) -- (-0.5,0);

\draw         (0.5,0) -- (2.5,0);
\draw[dashed] (3  ,0) -- (4  ,0);
\draw         (4.5,0) -- (5.5,0);
\draw         (0.5,2) -- (2.5,2);
\draw[dashed] (3  ,2) -- (4  ,2);
\draw         (4.5,2) -- (5.5,2);

\draw (1,-1) -- (1,1);
\draw (2,-1) -- (2,1);
\draw (5,-1) -- (5,1);
\draw (1, 3) -- (1,1);
\draw (2, 3) -- (2,1);
\draw (5, 3) -- (5,1);

\node at (6,1){$\bm{=}$};

\draw         (6.5,0) -- (7.5 ,0);
\draw[dashed] (8  ,0) -- (9   ,0);
\draw         (9.5,0) -- (10.5,0);
\draw        (10.5,0) -- (11.5,0);
\draw         (6.5,2) -- (7.5 ,2);
\draw[dashed] (8  ,2) -- (9   ,2);
\draw         (9.5,2) -- (10.5,2);
\draw        (10.5,2) -- (11.5,2);

\draw (7 ,-1) -- (7 ,1);
\draw (10,-1) -- (10,1);
\draw (11,-1) -- (11,1);
\draw (7,  3) -- (7 ,1);
\draw (10, 3) -- (10,1);
\draw (11, 3) -- (11,1);

\draw (11.5,0) -- (12.5,2);
\draw (11.5,2) -- (12.5,0);
\draw (12.5,0) -- (13,0);
\draw (12.5,2) -- (13,2);

\filldraw[red] (0  ,1) circle  (3pt);
\filldraw[red] (12 ,1) circle  (3pt);

\filldraw[blue] (1 ,0) circle  (3pt);
\filldraw[blue] (2 ,0) circle  (3pt);
\filldraw[blue] (5 ,0) circle  (3pt);
\filldraw[blue] (7 ,0) circle  (3pt);
\filldraw[blue] (10,0) circle  (3pt);
\filldraw[blue] (11,0) circle  (3pt);

\filldraw[blue] (1 ,2) circle  (3pt);
\filldraw[blue] (2 ,2) circle  (3pt);
\filldraw[blue] (5 ,2) circle  (3pt);
\filldraw[blue] (7 ,2) circle  (3pt);
\filldraw[blue] (10,2) circle  (3pt);
\filldraw[blue] (11,2) circle  (3pt);

\node[left] at (-0.5,1){$S_{ab}(u)$};
\node[right] at (12.5,1){$S_{ab}(u)$};
\node[below left] at (1,0){$M_{a}(u)$};
\node[above left] at (1,2){$M_{b}(u)$};
\node[below right] at (11,0){$M_{a}(u)$};
\node[above right] at (11,2){$M_{b}(u)$};

\node[below] at (1 ,-1){$1$};
\node[below] at (2 ,-1){$2$};
\node[below] at (5 ,-1){$L$};
\node[below] at (7 ,-1){$1$};
\node[below] at (10,-1){$L-1$};
\node[below] at (11,-1){$L$};

\node[left] at (-1,2){$a$};
\node[left] at (-1,0){$b$};
\node[right] at (13,2){$a$};
\node[right] at (13,0){$b$};

\end{tikzpicture} 
\end{center}

\caption{A graphical representation of the fundamental relation (\ref{FR}).
It corresponds to a representation of the Yang-Baxter equation (\ref{YBE-SR})
applied to an entire row of the lattice. The fundamental relation
provides the commutation relations between the monodromy matrix elements,
which are needed to compute the action of the transfer matrix on the
excited states. It also provides a sufficient condition for the transfer
matrix to commute with itself for different values of the spectral
parameter.}

\label{Fig.Mono-P}
\end{figure}
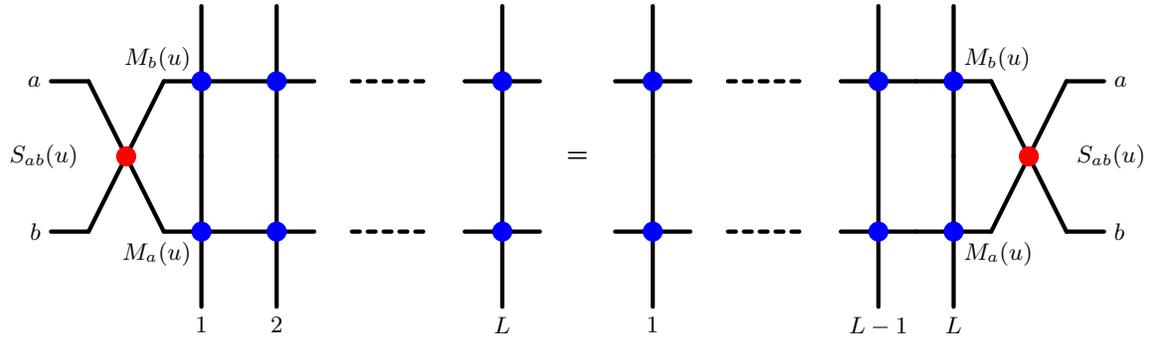

To evaluate this, the commutation relations between the diagonal operators
$A(u)$ and $D(u)$ with the creator operators $B(u_{k})$ $(1\leqslant k\leqslant n)$
must be known. All these commutation relations are provided by the
so-called \emph{fundamental relation} of the algebraic Bethe Ansatz,
\begin{equation}
S_{ab}(u-v)M_{a}(u)M_{b}(v)=M_{a}(v)M_{b}(u)S_{ab}(u-v).\label{FR}
\end{equation}
which can be thought as a representation of the Yang-Baxter equation
(\ref{YBE-SR}) applied to an entire row of the lattice \textendash{}
see Fig. \ref{Fig.Mono-P}. 

The fundamental relation (\ref{FR}) provides the commutation relation
between all the elements of the monodromy matrix (which must be written
in the Lax representation in order to the equation (\ref{FR}) make
sense). For the execution of the algebraic Bethe Ansatz, however,
only the following three relations are necessary: 
\begin{align}
A(u)B(v) & =a_{1}(u,v)B(v)A(u)+a_{2}(u,v)B(u)A(v),\label{AB6VP}\\
D(u)B(v) & =d_{1}(u,v)B(v)A(u)+d_{2}(u,v)B(u)D(v),\label{DB6VP}\\
B(u)B(v) & =B(v)B(u),\label{BB6VP}
\end{align}
where we introduced the following amplitudes:
\begin{align}
a_{1}(u,v) & =\frac{r_{1}\left(v-u\right)}{r_{2}\left(v-u\right)}, & a_{2}(u,v) & =-\frac{r_{3}\left(v-u\right)}{r_{2}\left(v-u\right)}, &  & \mathrm{and} & d_{1}(u,v) & =\frac{r_{1}\left(u-v\right)}{r_{2}\left(u-v\right)}, & d_{2}(u,v) & =-\frac{r_{3}\left(u-v\right)}{r_{2}\left(u-v\right)}.
\end{align}

The fundamental relation of the algebraic Bethe Ansatz (\ref{FR})
also reveals a profound property of models that have an $R$-matrix
satisfying the Yang-Baxter equation. To see that, take the trace of
(\ref{FR}) in both the $H_{a}$ and $V_{q}$ spaces. After we use
the properties of the trace and the definition of the transfer matrix
given by (\ref{tP}), we shall get directly the relation 
\begin{equation}
T(u)T(v)-T(v)T(u)=0,
\end{equation}
which shows us that the transfer matrix commutes with itself for different
values of the spectral parameter. This means that the transfer matrix
can be thought as the generator of infinitely many conserved quantities
in evolution \textendash{} the Hamiltonian being one of them \textendash{}
or, in other words, that the model is\emph{ integrable.}

\subsubsection{The computation of the eigenvalues and eigenstates of the transfer
matrix}

Now we have all that is needed to compute the action of the transfer
matrix on the $n$th excited state. This task, however, is not easy
at all. As a matter of a fact it is the most difficult step in the
execution of the algebraic Bethe Ansatz. Indeed, to compute this action
we need to use repeatedly the commutation relations (\ref{AB6VP})
and (\ref{DB6VP}), so that we can pass the diagonal operators $A$
and $D$ over all the creator operators $B$, after which we can act
with the diagonal operators on the reference state $\Psi_{0}$. It
is necessary, therefore, to use each commutation relation $n$ times,
which will generate $2^{n}$ terms each. 

Fortunately, after we analyze the first cases (e.g., for $n=1$ and
$n=2$), we might see that some patterns arises. In fact, we might
see that that some terms can be simplified and others can be grouped
together, so that the we can write a direct formula for the repeated
use of the commutation relations. These formulas are: 
\begin{align}
A(u)\prod_{k=1}^{n}B(u_{k}) & =\prod_{k=1}^{n}a_{1}(u,u_{k})B(u_{k})A(u)+B(u)\sum_{j=1}^{n}a_{2}(u,u_{k})\prod_{k=1,k\neq j}^{n}a_{1}(u_{j},u_{k})B(u_{k})A(u_{j}),\label{ABB-P}\\
D(u)\prod_{k=1}^{n}B(u_{k}) & =\prod_{k=1}^{n}d_{1}(u,u_{k})B(u_{k})D(u)+B(u)\sum_{j=1}^{n}d_{2}(u,u_{k})\prod_{k=1,k\neq j}^{n}d_{1}(u_{j},u_{k})B(u_{k})D(u_{j}).\label{DBB-P}
\end{align}
A rigorous proof of (\ref{ABB-P}) and (\ref{DBB-P}) is, nevertheless,
required. This proof follows from mathematical induction: the assumption
is clearly true for $n=1$. Assume that it also holds for general
$n$. Then, for $n+1$ we have that,
\begin{align}
A(u)\prod_{k=1}^{n+1}B(u_{k}) & =a_{1}(u,u_{n+1})B(u_{n+1})A(u)\prod_{k=1}^{n}B(u_{k})-a_{2}(u,u_{n+1})B(u)A(u_{n+1})\prod_{k=1}^{n}B(u_{k}),\\
D(u)\prod_{k=1}^{n+1}B(u_{k}) & =d_{1}(u,u_{n+1})B(u_{n+1})D(u)\prod_{k=1}^{n}B(u_{k})-d_{2}(u,u_{n+1})B(u)D(u_{n+1})\prod_{k=1}^{n}B(u_{k}),
\end{align}
 where we made use of the fact the $B$ operators commute with themselves
in order to put $B(u_{n+1})$ on the left. Now, using (\ref{ABB-P})
we get that 
\begin{align}
A(u)\prod_{k=1}^{n+1}B(u_{k}) & =\prod_{k=1}^{n+1}a_{1}(u,u_{k})B(u_{k})A(u)\nonumber \\
 & +a_{1}(u,u_{n+1})B(u_{n+1})\sum_{j=1}^{n}a_{2}(u,u_{j})B(u)\prod_{k=1,k\neq j}^{n}a_{1}(u_{j},u_{k})B(u_{k})A(u_{j})\nonumber \\
 & +a_{2}(u,u_{n+1})B(u)\prod_{k=1}^{n}a_{1}(u_{n+1},u_{k})B(u_{k})A(u_{n+1})\nonumber \\
 & +a_{2}(u,u_{n+1})B(u)\sum_{j=1}^{n}a_{2}(u_{n+1},u_{j})B(u_{n+1})\prod_{k=1,k\neq j}^{n}a_{1}(u_{j},u_{k})B(u_{k})A(u_{j}),
\end{align}
 and, similarly, using (\ref{DBB-P}), we get, 
\begin{align}
D(u)\prod_{k=1}^{n+1}B(u_{k}) & =\prod_{k=1}^{n+1}d_{1}(u,u_{k})B(u_{k})D(u)\nonumber \\
 & +d_{1}(u,u_{n+1})B(u_{n+1})\sum_{j=1}^{n}d_{2}(u,u_{j})B(u)\prod_{k=1,k\neq j}^{n}d_{1}(u_{j},u_{k})B(u_{k})D(u_{j})\nonumber \\
 & +d_{2}(u,u_{n+1})B(u)\prod_{k=1}^{n}d_{1}(u_{n+1},u_{k})B(u_{k})D(u_{n+1})\nonumber \\
 & +d_{2}(u,u_{n+1})B(u)\sum_{j=1}^{n}d_{2}(u_{n+1},u_{j})B(u_{n+1})\prod_{k=1,k\neq j}^{n}d_{1}(u_{j},u_{k})B(u_{k})D(u_{j}).
\end{align}
Then, we can see that the second and the fourth terms in each expressions
above can be grouped together and, after we use the identities, 
\begin{align}
a_{1}(u,u_{n+1})a_{2}(u,u_{j})+a_{2}(u,u_{n+1})a_{2}(u_{n+1},u_{j}) & =a_{2}(u,u_{j})a_{1}(u_{j},u_{n+1}),\qquad1\leqslant j\leqslant n,\label{ID6P1}\\
d_{1}(u,u_{n+1})d_{2}(u,u_{j})+d_{2}(u,u_{n+1})d_{2}(u_{n+1},u_{j}) & =d_{2}(u,u_{j})d_{1}(u_{j},u_{n+1}),\qquad1\leqslant j\leqslant n,\label{ID6P2}
\end{align}
those terms can also be grouped with the third ones in (\ref{ABB-P})
and (\ref{DBB-P}), after we extend the summation to $n+1$. This
lead us to same expressions (\ref{ABB-P}) and (\ref{DBB-P}) with
$n+1$ in the place of $n$, which proves the assumption.

Now, from (\ref{ABB-P}) and (\ref{DBB-P}) it is an easy matter to
compute the action of the transfer matrix on the $n$th excited state.
It happens that the action of $T(u)$ on $\Psi_{n}\left(u_{1},\ldots,u_{n}\right)$
can be written as,
\begin{equation}
T(u)\Psi_{n}\left(u_{1},\ldots,u_{n}\right)=\tau_{n}\left(u|u_{1},\ldots,u_{n}\right)\Psi_{n}\left(u_{1},\ldots,u_{n}\right)+\sum_{k=1}^{n}\beta_{n}^{k}\left(u|u_{1},\ldots,u_{n}\right)\Psi_{n}\left(u_{k}^{\times}\right),\label{ActionTP}
\end{equation}
 where,
\begin{equation}
\tau_{n}\left(u|u_{1},\ldots,u_{n}\right)=\alpha(u)\prod_{k=1}^{n}a_{1}\left(u,u_{k}\right)+\delta(u)\prod_{k=1}^{n}d_{1}\left(u,u_{k}\right),\label{eigenP}
\end{equation}
\begin{align}
\beta_{n}^{k}\left(u|u_{1},\ldots,u_{n}\right) & =\alpha\left(u_{k}\right)a_{2}\left(u,u_{k}\right)\prod_{i=1,i\neq k}^{n}a_{1}\left(u_{k},u_{i}\right)+\delta\left(u_{k}\right)d_{2}\left(u,u_{k}\right)\prod_{i=1,i\neq k}^{n}d_{1}\left(u_{k},u_{i}\right), &  & 1\leqslant k\leqslant n,
\end{align}
and we introduced the notation 
\begin{equation}
\Psi_{n}\left(u_{k}^{\times}\right)=B(u)\prod_{j=1,j\neq k}^{n}B(u_{j})\Psi_{0}.
\end{equation}

The main point here is that the first term in (\ref{ActionTP}) is
proportional to $\Psi_{n}\left(u_{1},\ldots,u_{n}\right)$, while
the other terms are not. Therefore, the requirement that the transfer
matrix satisfies an eigenvalue equation means that all the \emph{unwanted
terms} \textendash{} i.e., the terms which are not proportional to
$\Psi_{n}\left(u_{1},\ldots,u_{n}\right)$ \textendash{} need to vanish.
This can be ensured by fixing appropriated values for the rapidities
$u_{1},\ldots,u_{k}$, which until now were arbitrary. In fact, as
we impose that the following system of non-linear equations, 
\begin{align}
\beta_{n}^{k}\left(u|u_{1},\ldots,u_{n}\right) & =\alpha\left(u_{k}\right)a_{2}\left(u,u_{k}\right)\prod_{i=1,i\neq k}^{n}a_{1}\left(u_{k},u_{i}\right)+\delta\left(u_{k}\right)d_{2}\left(u,u_{k}\right)\prod_{i=1,i\neq k}^{n}d_{1}\left(u_{k},u_{i}\right)=0, &  & 1\leqslant k\leqslant n,\label{BAEP}
\end{align}
is satisfied, all the unwanted terms will vanish. The system of equations
are called the \emph{Bethe Ansatz equations} of the six-vertex model
with periodic boundary conditions \cite{Bethe1931}. Their solutions
provide the correct values for the rapidities $u_{1},\ldots u_{n}$
necessary to vanish the unwanted terms. 

\subsubsection{The solutions of the Bethe Ansatz equations}

We remark that if were possible to solve the Bethe Ansatz equations
analytically, then we would obtain a complete analytical answer for
the spectral problem considered by the Bethe Ansatz technique but,
unfortunately, they are too complex for such an ambitious endeavor
be accomplished by now \textendash{} it is only up to the second excited
state that a complete analytical solution of the Bethe Ansatz equations
were obtained so far \cite{VieiraLimaSantos2015,Vieira2017} \textendash ,
reason by which the Bethe Ansatz equations are usually solved through
numerical methods \cite{BabelonVegaViallet1983,Baxter2002}.

\subsection{The combinatorial approach\label{Section-Combinatorial-P}}

In the previous section we discussed how the periodic algebraic Bethe
Ansatz is usually implemented. We have seen that it depends on several
steps, which might make the method seem difficult, especially for
introductory audiences. When some subject has several technical and
complicated details, it is always desirable, when possible, to set
up symbolic approach for it, so that its main features can be qualitatively
understood and the final results obtained in an easier way \textendash{}
a classical example of this is the use of Feynman's diagrams in quantum
field theory. In this section we shall describe such a symbolic method
for the algebraic Bethe Ansatz in terms of combinatorial diagrams
(for an analysis of the Bethe Ansatz through the tools of tensor-networks,
see \cite{MurgKorepinVerstraete2012}). Remember that the most laborious
step in the execution of the algebraic Bethe Ansatz is the computation
of the action of the transfer matrix on the $n$th excited state.
This computation relies on the repeated use of the commutation relations
(\ref{AB6VP}) and (\ref{DB6VP}). Or method consists in representing
these commutation relations by simple combinatorial diagrams, so that
the repeated use of them gives place to simple \emph{combinatorial
trees}. As we shall see in the sequel, the analysis of these combinatorial
trees provides in a straightforward way the eigenvalues and eigenstates
of the transfer matrix as well as the respective Bethe Ansatz equations. 

Before, however, we discuss how this combinatorial approach works,
we should present some definitions and nomenclatures that we shall
make extensive use in the sequel. These nomenclatures appear in the
theory of graphs, more specifically in the study of \emph{combinatorial
rooted trees}. Following \cite{West2001} and \cite{Bona2011}, a
\emph{rooted tree} is defined as a directed graph in which any two
\emph{nodes} (the vertices of the graph) are connected by exactly
one \emph{path}. A \emph{labeled tree} is a tree whose nodes are specified
by means of labels. The first node of the tree (which connects all
other nodes) is called its \emph{root} and the last node of a given
path of the tree is called a \emph{leaf} node. We say that a given
node is at the level $k$ of the tree if the path connecting this
node to the root contains exactly $k+1$ nodes (counting both the
root and the referred node, so that the root is always at the level
zero). The \emph{length} of a path is defined as the number of nodes
it contains, from the root to the leaves, so that the length of the
tree is also the length of its longest path. If all paths of the tree
have the same length (as is the case for the trees considered here),
we shall refer to it as a \emph{pruned tree}. Moreover, we establish
a parental relationship between the nodes of the tree: given a node
at the level $k$ of the tree and another node at the level $k+1$,
we say that the first is the \emph{parent} of the second \textendash{}
and, accordingly, that the second is a \emph{child} of the first \textendash{}
if they belong to the same path of the tree. Finally, a disjoint union
of trees is usually called a \emph{forest}. 

The combinatorial approach we are going to present is actually very
simple: it consists in representing the commutation relations 
\begin{align}
A(u)B(v) & =a_{1}(u,v)B(v)A(u)+a_{2}(u,v)B(u)A(v),\label{CRABP}\\
D(u)B(v) & =d_{1}(u,v)B(v)A(u)+d_{2}(u,v)B(u)D(v),\label{CRDBP}
\end{align}
by simple diagrams with two outputs, as below: \begin{equation}
\begin{tikzpicture}[edge from parent fork down,
level 1/.style={sibling distance=1cm},
level distance=0.7cm
]
\tikzstyle{hollow circle}=[circle,    thick, draw, inner sep=2.5] 
\tikzstyle{solid  circle}=[circle,    thick, draw, inner sep=2.5, fill=blue]
\tikzstyle{hollow square}=[rectangle, thick, draw, inner sep=3] 
\tikzstyle{solid  square}=[rectangle, thick, draw, inner sep=3, fill=blue]
\tikzstyle{edge from parent}=[thick,draw]
\node(0)[hollow square,red]{\textsc{a}} 
  child {node(1)[hollow circle,blue]{} }
  child {node(2)[solid  circle,blue]{} }
; 
\node[right] at (1) {$\hspace{0,1cm} u$};
\node[right] at (2) {$\hspace{0,1cm} v$};
\end{tikzpicture}
\hspace{1cm} \text{and} \hspace{1cm}
\begin{tikzpicture}[edge from parent fork down,
level 1/.style={sibling distance=1cm},
level distance=0.7cm
]
\tikzstyle{hollow circle}=[circle,    thick, draw, inner sep=2.5] 
\tikzstyle{solid  circle}=[circle,    thick, draw, inner sep=2.5, fill=blue]
\tikzstyle{hollow square}=[rectangle, thick, draw, inner sep=3] 
\tikzstyle{solid  square}=[rectangle, thick, draw, inner sep=3, fill=blue]
\tikzstyle{edge from parent}=[thick,draw]
\node(0)[hollow square, red]{\textsc{d}} 
  child {node(1)[hollow circle,blue]{} }
  child {node(2)[solid  circle,blue]{} }
; 
\node[right] at (1) {$\hspace{0,1cm} u$};
\node[right] at (2) {$\hspace{0,1cm} v$};
\end{tikzpicture}
\label{diagram1P}
\end{equation} Notice that in these diagrams the root indicates what commutation
relation we are talking about. In each diagram, the nodes (hollow
or filled) represent one of the two terms in the respective commutation
relations. For instance, we can assume that the hollow node represents
the first term in the commutation relations (\ref{CRABP}) or (\ref{CRDBP}),
while the filled node represents the second term in (\ref{CRABP})
or (\ref{CRDBP}). Thus, we can write something like this: \begin{equation}
\begin{tabular}{c}
\begin{tikzpicture}
\node(A) at (0,0) [rectangle, red, thick, draw, inner sep=2.5]{\textsc{a}}; 
\node(1) at (1,0) [circle, blue,   thick, draw, inner sep=2.5]{};
\node[right] at (1){$\hspace{0.1cm}\equiv a_{1}(u,v)B(v)A(u),$};
\draw[thick,cap=round] (A) -- (1);
\end{tikzpicture}
\qquad
\begin{tikzpicture}
\node(A) at (0,0) [rectangle, red, thick, draw, inner sep=2.5]{\textsc{d}}; 
\node(1) at (1,0) [circle, blue,   thick, draw, inner sep=2.5]{};
\node[right] at (1){$\hspace{0.1cm}\equiv d_{1}(u,v)B(v)D(u),$};
\draw[thick,cap=round] (A) -- (1);
\end{tikzpicture}\tabularnewline
\begin{tikzpicture}
\node(A) at (0,0) [rectangle, red, thick, draw, inner sep=2.5]{\textsc{a}}; 
\node(1) at (1,0) [circle, blue, fill=blue,  thick, draw, inner sep=2.5]{};
\node[right] at (1){$\hspace{0.1cm}\equiv a_{2}(u,v)B(u)A(v),$};
\draw[thick,cap=round] (A) -- (1);
\end{tikzpicture}
\qquad
\begin{tikzpicture}
\node(A) at (0,0) [rectangle, red, thick, draw, inner sep=2.5]{\textsc{d}}; 
\node(1) at (1,0) [circle, blue, fill=blue,  thick, draw, inner sep=2.5]{};
\node[right] at (1){$\hspace{0.1cm}\equiv d_{2}(u,v)B(u)D(v).$};
\draw[thick,cap=round] (A) -- (1);
\end{tikzpicture}\tabularnewline
\end{tabular}\end{equation} Finally, the labels $u$ and $v$ attached at the
side of each node in the diagrams (\ref{diagram1P}) indicate the
argument of the diagonal operators, $A$ or $D$, in the respective
term of the commutation relation (we shall see that they are very
useful in the general case). 

Now, the repeated use of the commutation relations can be represented
by a \emph{labeled combinatorial tree}. For example, let us consider
the computation of $A(u_{0})B(u_{1})B(u_{2})\cdots B(u_{n})$, which
appears in the action of the operator $A(u_{0})$ on the $n$th excited
state $\Psi_{n}(u_{1},\ldots,u_{n})$. This can be represented\footnote{In this section and in section \ref{Section-Combinatorial-B}, we
shall usually write the spectral parameter $u$ as $u_{0}$.} by the following pruned binary tree of length $n$: \begin{equation}
\begin{tikzpicture}[edge from parent fork down,
level 1/.style={sibling distance=4cm}, 
level 2/.style={sibling distance=2cm},
level 3/.style={sibling distance=1cm},
level 4/.style={sibling distance=0.5cm},
level distance=0.7cm
]
\tikzstyle{hollow circle}=[circle, blue,   thick, draw, inner sep=2.5] 
\tikzstyle{solid  circle}=[circle, blue,   thick, draw, inner sep=2.5, fill=blue]
\tikzstyle{hollow square}=[rectangle,red, thick, draw, inner sep=3] 
\tikzstyle{solid  square}=[rectangle,red, thick, draw, inner sep=3, fill=blue]
\tikzstyle{edge from parent}=[thick,draw]

\node(0)[hollow square]{\textsc{a}}
  child{node(1)[hollow circle]{}
    child{node(11)[hollow circle]{}
      child{node(111)[hollow circle]{}}
      child{node(112)[solid  circle]{}}
    }
    child{node(12)[solid  circle]{}
      child{node(121)[hollow circle]{}}
      child{node(122)[solid circle]{}}
    }
  } 
 child{node(2)[solid circle]{}
    child{node(21)[hollow circle]{}
      child{node(211)[hollow circle]{}}
      child{node(212)[solid  circle]{}}
    }
    child{node(22)[solid  circle]{}
      child{node(221)[hollow circle]{}}
      child{node(222)[solid circle]{}}
    }
  } 
; 

\node[right] at  (1) {$\hspace{0,1cm} u_0$};
\node[right] at  (2) {$\hspace{0,1cm} u_1$};
\node[right] at (11) {$\hspace{0,1cm} u_0$};
\node[right] at (12) {$\hspace{0,1cm} u_2$};
\node[right] at (21) {$\hspace{0,1cm} u_1$};
\node[right] at (22) {$\hspace{0,1cm} u_2$};
\node[right] at (111) {$\hspace{0,1cm} u_0$};
\node[right] at (112) {$\hspace{0,1cm} u_3$};
\node[right] at (121) {$\hspace{0,1cm} u_2$};
\node[right] at (122) {$\hspace{0,1cm} u_3$};
\node[right] at (211) {$\hspace{0,1cm} u_1$};
\node[right] at (212) {$\hspace{0,1cm} u_3$};
\node[right] at (221) {$\hspace{0,1cm} u_2$};
\node[right] at (222) {$\hspace{0,1cm} u_3$};

\node[below] at (111) {$\vdots$};
\node[below] at (112) {$\vdots$};
\node[below] at (121) {$\vdots$};
\node[below] at (122) {$\vdots$};
\node[below] at (211) {$\vdots$};
\node[below] at (212) {$\vdots$};
\node[below] at (221) {$\vdots$};
\node[below] at (222) {$\vdots$};
\end{tikzpicture}
\label{diagramNP}
\end{equation}and an identical diagram holds for the computation of $D(u_{0})B(u_{1})B(u_{2})\cdots B(u_{n})$.

As before, the hollow nodes in the diagram (\ref{diagramNP}) mean
that the first term of the commutation relations (\ref{CRABP}) or
(\ref{CRDBP}) is gathered at that point, while the filled nodes mean
that the second term is kept there. The diagram will contain, therefore,
$2^{n}$ paths, each one of them representing a term that would be
obtained from the algebraic Bethe Ansatz. Any path $P$ of the diagram
can be specified by the set $\eta=\left(\eta_{1},\ldots,\eta_{n}\right)$,
where the parameters $\eta_{k}$ $(1\leqslant k\leqslant n)$ can
assume only the values $1$ or $2$, depending on whether the node
in the level $k$ of the path $P$ is hollow $(\eta_{k}=1)$ or filled
$(\eta_{k}=2)$. 

The label $\lambda$ at the side of every node represents the argument
of the diagonal operators, $A$ or $D$, on that point of the respective
diagram. These labels can be easily obtained following the simple
rule: 
\begin{itemize}
\item Hollow nodes always inherit the label of his parent, while the label
$\lambda_{k}^{P}$ of any filled node on the level $k$ of the path
$P$ is always $u_{k}$ (the label of the root being defined as $u_{0}$). 
\end{itemize}
This rule comes from a direct analysis of the commutation relations
(\ref{CRABP}) and (\ref{CRDBP}). In fact, starting, for instance,
with $A(u_{0})B(u_{1})B(u_{2})\cdots B(u_{n})$, we may notice that
the use of the commutation relation (\ref{CRABP}) has always the
effect of permuting the operators $A$ and $B$. If we gather, in
the one hand, the first term of the commutation relation (\ref{CRABP}),
then we shall get a quantity proportional to $B(u_{1})A(u_{0})B(u_{2})\cdots B(u_{n})$,
where the argument of the diagonal operator $A$ is still the same
as before. If we gather, on the other hand, the second term in (\ref{CRABP}),
then we would get a quantity proportional to $B(u_{0})A(u_{1})B(u_{2})\cdots B(u_{n})$,
so that the argument of the diagonal operator $A$ is permuted with
that of the operator $B$. We can say, therefore, that a hollow node
does not change the argument of the diagonal operators $A$ or $D$,
so that their labels must be the same as the labels of their parent,
while a filled node permutes the argument of the diagonal operators
$A$ or $D$ with that of the creator operator $B$ in that point,
that is, the label of any filled node at the level $k$ of the diagram
must be $u_{k}$. This means that the label $\lambda_{k}^{P}$ of
a node at the level $k$ of a given path $P$ in any diagram ($A$
or $D$) can be determined recursively by the formula 
\begin{equation}
\lambda_{k}^{P}=\begin{cases}
\lambda_{k-1}^{P}, & \eta_{k}=1,\\
u_{k}, & \eta_{k}=2,
\end{cases}
\end{equation}
 provided we define the label of the root as $\lambda_{0}^{P}=u_{0}$.

The key point of this approach is that there is a one-to-one correspondence
between each path of the diagrams and each one of the final terms
obtained from the usual algebraic Bethe Ansatz\footnote{As an example, we present below the paths associated with the $A$
diagram for the second excited state ($n=2$) and the corresponding
mathematical expressions: \begin{center} %
\begin{tabular}{l}
\begin{tikzpicture}

\node(A) at (0,0) [rectangle, red, thick, draw, inner sep=2.5]{\textsc{a}}; 
\node(B) at (1,0) [circle,   blue, thick, draw, inner sep=2.5]{};
\node(C) at (2,0) [circle,   blue, thick, draw, inner sep=2.5]{};
\draw[thick] (A) -- (B);
\draw[thick] (B) -- (C);
\node[right] at (C){\hspace{0.1cm}$ \equiv a_{1}(u_{0},u_{1})a_{1}(u_{0},u_{2})\alpha\left(u_{0}\right)B(u_{1})B(u_{2})\Psi_{0},$}; 
\end{tikzpicture}
\qquad
\begin{tikzpicture}
\node(A) at (0,0) [rectangle, red, thick, draw, inner sep=2.5]{\textsc{a}}; 
\node(B) at (1,0) [circle,   blue, thick, draw, inner sep=2.5]{};
\node(C) at (2,0) [circle,   blue, fill=blue,thick, draw, inner sep=2.5]{};
\draw[thick] (A) -- (B);
\draw[thick] (B) -- (C);
\node[right] at (C){\hspace{0.1cm}$ \equiv a_{1}(u_{0},u_{1})a_{2}(u_{0},u_{2})\alpha\left(u_{2}\right)B(u_{0})B(u_{1})\Psi_{0},$};  
\end{tikzpicture}\tabularnewline
\begin{tikzpicture}
\node(A) at (0,0) [rectangle, red, thick, draw, inner sep=2.5]{\textsc{a}}; 
\node(B) at (1,0) [circle,   blue, fill=blue,thick, draw, inner sep=2.5]{};
\node(C) at (2,0) [circle,   blue, thick, draw, inner sep=2.5]{};
\draw[thick] (A) -- (B);
\draw[thick] (B) -- (C);
\node[right] at (C){\hspace{0.1cm}$ \equiv a_{2}(u_{0},u_{1})a_{1}(u_{1},u_{2})\alpha\left(u_{1}\right)B(u_{0})B(u_{2})\Psi_{0},$};  
\end{tikzpicture}

\qquad

\begin{tikzpicture}
\node(A) at (0,0) [rectangle, red, thick, draw, inner sep=2.5]{\textsc{a}}; 
\node(B) at (1,0) [circle,   blue, fill=blue, thick, draw, inner sep=2.5]{};
\node(C) at (2,0) [circle,   blue, fill=blue, thick, draw, inner sep=2.5]{};
\draw[thick] (A) -- (B);
\draw[thick] (B) -- (C);
\node[right] at (C){\hspace{0.1cm}$ \equiv a_{2}(u_{0},u_{1})a_{2}(u_{1},u_{2})\alpha\left(u_{2}\right)B(u_{0})B(u_{1})\Psi_{0}.$}; 
\end{tikzpicture}\tabularnewline
\end{tabular}\end{center} The paths associated with the $D$ diagram are the same
as the above ones, except that the coefficients $a_{1}$ and $a_{2}$
should be replaced by $d_{1}$ and $d_{2}$, respectively.}. In fact, after we use one of the commutation relations (\ref{CRABP})
or (\ref{CRDBP}) $n$ times, we shall obtain $2^{n}$ terms, or \emph{values}
$V$, each of which is of the form 
\begin{equation}
V=W\left|S\right\rangle ,
\end{equation}
where $W$ (the \emph{weight} of the term $T$, or the path $P$ it
represents) corresponds to a specific product of the coefficients
appearing on the commutation relations, times the action of the diagonal
operator on the reference state $\Psi_{0}$, and $\left|S\right\rangle $
(the \emph{state} associated with the term $T$ or the corresponding
path $P$) consists in a given product of the $B$ operators. The
weight of a path can be found defining the contribution of each node
it contains, plus a contribution of the leaf node. In this way, we
define the weights of a given node and the contribution of the leaf
node as follows: 
\begin{itemize}
\item The weight of a \emph{hollow} node in the level $k$ of a given path
$P_{A}$ $\left[P_{D}\right]$ equals $a_{1}(\lambda_{k-1}^{P},u_{k})$
$\left[d_{1}(\lambda_{k-1}^{P},u_{k})\right]$, while the weight of
a given \emph{filled} node in the level $k$ of a given path $P_{A}$
$\left[P_{D}\right]$ equals $a_{2}(\lambda_{k-1}^{P},u_{k})$ $\left[d_{2}(\lambda_{k-1}^{P},u_{k})\right]$. 
\item The leaf node of the path $P_{A}$ $\left[P_{D}\right]$ contributes
to the weight of the respective path with the factor $\alpha(\lambda_{n}^{P})$
$\left[\delta(\lambda_{n}^{P})\right]$, which arises from the action
of the operator $A(\lambda_{n})$ $\left[D(\lambda_{n})\right]$ on
$\Psi_{0}$.
\end{itemize}
In a more condensed way, the weight associated with each $P_{A}(\eta_{1},\ldots,\eta_{n})$
of the $A$ diagram, and the weight of each path $P_{D}(\eta_{1},\ldots,\eta_{n})$
of the $D$ diagram, are given by,
\begin{align}
W\left(P_{A}(\eta_{1},\ldots,\eta_{n})\right) & =\alpha(\lambda_{n}^{P})\prod_{k=1}^{n}a_{\eta_{k}}\left(\lambda_{k-1}^{P},u_{k}\right), & W\left(P_{D}(\eta_{1},\ldots,\eta_{n})\right) & =\delta(\lambda_{n}^{P})\prod_{k=1}^{n}d_{\eta_{k}}\left(\lambda_{k-1}^{P},u_{k}\right).\label{WPA}
\end{align}

Similarly, the state associated with each path can be determined by
a single rule: 
\begin{itemize}
\item The state $|P^{(u_{k})}\rangle$ associated with any path $P$ whose
leaf node has the label $\lambda_{n}^{P}=u_{k}$ is given by the product
of all operators $B(u_{j})$, $0\leqslant j\neq k\leqslant n$, times
$\Psi_{0}$. That is, 
\begin{equation}
|P^{\left(u_{k}\right)}\rangle=\prod_{j=0,j\neq k}^{n}B(u_{j})\Psi_{0}.\label{Puk}
\end{equation}
\end{itemize}
In fact, the condition for a given path to end with the label $u_{k}$
is that it contains a filled node in the level $k$ and no other filled
node in the higher levels \textendash{} which is a consequence of
the rule determining the labels of the diagrams. Therefore, in both
the $A$ and $D$ diagrams there will be only one path ending with
the label $u_{0}$ and exactly $2^{k-1}$ paths that end with the
label $u_{k}$ for $k\geqslant1$. 

From what was said above, it is an easy matter to determine the whole
action of the transfer matrix on the $n$th excited state. In fact,
the action of the $T(u_{0})$ on $\Psi_{n}\left(u_{1},\ldots,u_{n}\right)$
will be given by the sum of the values $V(P)=W(P)\left|P\right\rangle $
associated with all paths $P$ of the $A$ and $D$ diagrams (i.e.,
by a sum over the \emph{forest}, in the jargon of graph theory). That
is, 
\begin{equation}
T(u_{0})\Psi_{n}\left(u_{1},\ldots,u_{n}\right)=\sum_{\eta_{1},\ldots,\eta_{n}=1}^{2}\left[V_{A}\left(\eta_{1},\ldots,\eta_{n}\right)+V_{D}\left(\eta_{1},\ldots,\eta_{n}\right)\right]\left|P\left(\eta_{1},\ldots,\eta_{n}\right)\right\rangle .\label{ActionTPcomb}
\end{equation}
Gathering all the paths with the same state (that is, collecting all
the paths terminating with the same label), we can also write this
as 
\begin{align}
T(u_{0})\Psi_{n}\left(u_{1},\ldots,u_{n}\right) & =\left[W\left(P_{A}^{(u_{0})}\right)+W\left(P_{D}^{(u_{0})}\right)\right]\left|P^{(u_{0})}\right\rangle +\sum_{k=1}^{n}\left[W\left(P_{A}^{(u_{k})}\right)+W\left(P_{D}^{(u_{k})}\right)\right]\left|P^{(u_{k})}\right\rangle .\label{ActionWP}
\end{align}
This corresponds to partitioning the $2^{n}$ terms in expression
(\ref{ActionTPcomb}) into the $n+1$ terms in expression (\ref{ActionWP}).
Notice that this partition is only possible thanks to the identity
$1+\left(1+2+\cdots+2^{n-1}\right)=2^{n}$. Identifying
\begin{equation}
\tau_{n}\left(u_{0}|u_{1},\ldots,u_{n}\right)\Psi_{n}\left(u_{1},\ldots,u_{n}\right)\equiv\left[W\left(P_{A}^{(u_{0})}\right)+W\left(P_{D}^{(u_{0})}\right)\right]\left|P^{(u_{0})}\right\rangle ,
\end{equation}
\begin{equation}
\beta_{n}^{k}\left(u_{0}|u_{1},\ldots,u_{n}\right)\Psi_{n}\left(u_{k}^{\times}\right)\equiv\left[W\left(P_{A}^{(u_{k})}\right)+W\left(P_{D}^{(u_{k})}\right)\right]\left|P^{(u_{k})}\right\rangle ,\qquad1\leqslant k\leqslant n,
\end{equation}
 and 
\begin{equation}
\left|P^{(u_{k})}\right\rangle \equiv\Psi_{n}\left(u_{k}^{\times}\right)=\prod_{j=0,j\neq k}^{n}B(u_{j}),\qquad1\leqslant k\leqslant n,
\end{equation}
 we can rewrite (\ref{ActionWP}) as 
\begin{equation}
T(u_{0})\Psi_{n}\left(u_{1},\ldots,u_{n}\right)=\tau_{n}\left(u_{0}|u_{1},\ldots,u_{n}\right)\Psi_{n}\left(u_{1},\ldots,u_{n}\right)+\sum_{k=1}^{n}\beta_{n}^{k}\left(u_{0}|u_{1},\ldots,u_{n}\right)\Psi_{n}\left(u_{k}^{\times}\right),
\end{equation}
which provides a complete agreement with (\ref{ActionTP}) for the
action of the transfer matrix on the $n$th excited state.

The analysis above gives a complete description for the action of
the transfer matrix $T(u_{0})$ on the $n$th excited state $\Psi_{n}$.
In practice, however, we are usually interested only in the eigenvalues
of the transfer matrix and in the Bethe Ansatz equations. Let us show
now how the eigenvalues can be obtained in a straightforward way through
the analysis of the diagrams. 

The eigenvalues of the transfer matrix can be determined through the
analysis of the diagrams by following only one simple rule: 
\begin{itemize}
\item The eigenvalues of the transfer matrix are determined by the sum of
the weights over all paths of the diagrams $A$ and $D$ that end
with the label $u_{0}$. 
\end{itemize}
In fact, only in this case the state of the path will be proportional
to $\Psi_{n}(u_{1},\ldots,u_{n})$. Notice that a path will end with
the label $u_{0}$ if, and only if, it contains no filled nodes \textendash{}
i.e., if it contains only hollow nodes. Since, however, there is only
one path in each diagram satisfying this requirement \textendash{}
namely, the paths $P_{A}\left(1,\ldots,1\right)$ and $P_{D}\left(1,\ldots,1\right)$
\textendash , we conclude that the eigenvalues are determined by \begin{equation}
\begin{tikzpicture}
\node(T) at (-3.05,0){$\tau_n\left(u_0|u_1,\ldots,u_n\right)\Psi_n\left(u_1,\ldots,u_n\right)=$};
\node(A) at (0,0) [rectangle, red, thick, draw, inner sep=2.5]{\textsc{a}};
\node(B) at (1,0) [circle, blue  , thick, draw, inner sep=2.5]{};
\node(C) at (2,0) [circle, blue   , thick, draw, inner sep=2.5]{};

\draw[thick] (A) -- (B);
\draw[thick] (B) -- (C);

\node(P) at (2.45,0){$+$};
\node(E) at (3,0) [rectangle,red, thick, draw, inner sep=2.5]{\textsc{d}};
\node(F) at (4,0) [circle, blue   , thick, draw, inner sep=2.5]{};
\node(G) at (5,0) [circle, blue   , thick, draw, inner sep=2.5]{};

\draw[thick] (E) -- (F);
\draw[thick] (F) -- (G);

\end{tikzpicture}
\end{equation}which means, according to the weights of these paths, that
\begin{equation}
\tau_{n}\left(u_{0}|u_{1},\ldots,u_{n}\right)=\alpha(u_{0})\prod_{k=1}^{n}a_{1}\left(u_{0},u_{k}\right)+\delta(u_{0})\prod_{k=1}^{n}d_{1}\left(u_{0},u_{k}\right).\label{EigenPnew}
\end{equation}

Similarly, the Bethe Ansatz equations can be determined by the rule:
\begin{itemize}
\item The Bethe Ansatz equation fixing the rapidity $u_{k}$ is determined
by a sum over all paths of the $A$ and $D$ diagrams whose leaf node
has the label $\lambda_{n}^{P}=u_{k}$.
\end{itemize}
In fact, we have seen that the state associated with a given path
ending with the label $u_{k}$ is just $\left|P^{\left(u_{k}\right)}\right\rangle =\Psi\left(u_{k}^{\times}\right)$.
The condition for a given path to end with the label $u_{k}$ is that
it contains a filled node at the level $k$ and no other filled node
at higher levels. In the lower levels, however, the nodes can be of
any type, which means that we should sum over all the possible types
of nodes in the lower levels of the diagrams. We can graphically express
this as follows: \begin{equation}
\begin{tikzpicture}[scale=1]
\node(t) at (0,0)[left]{$\beta_n^k\left(u_0|u_1,\ldots,u_n\right)\Psi_n\left(u_1,\ldots,u_n|u_{k}^{\times}\right)=\hspace{0.2cm}$};
\node(0) at (0,0) [rectangle, red, thick, draw, inner sep=2.5]{\textsc{a}};
\node(1) at (1,0) [circle, blue   ,fill=black!30, thick, draw, inner sep=2.5]{};
\node(2) at (2,0) [circle, blue   ,fill=black!30, thick, draw, inner sep=2.5]{};
\node(3) at (3,0) [circle, blue   ,fill=blue, thick, draw, inner sep=2.5]{};
\node(4) at (4,0) [circle, blue   , thick, draw, inner sep=2.5]{};
\node(5) at (5,0) [circle, blue   , thick, draw, inner sep=2.5]{};

\draw[thick] (0) -- (1);
\draw[thick] (1) -- (2);
\draw[thick,dashed] (2) -- (3);
\draw[thick] (3) -- (4);
\draw[thick] (4) -- (5);

\node(P)  at (5.45,0){$+$};
\node(ja)[below right] at(3,0){$k$};
\node(jd)[below right] at(9,0){$k$};
\node(6)  at (6,0) [rectangle, red, thick, draw, inner sep=2.5]{\textsc{d}};
\node(7)  at (7,0) [circle, blue   , fill=black!30, thick, draw, inner sep=2.5]{};
\node(8)  at (8,0) [circle, blue   , fill=black!30,thick, draw, inner sep=2.5]{};
\node(9)  at (9,0) [circle, blue   , fill=blue, thick, draw, inner sep=2.5]{};
\node(10) at (10,0) [circle, blue   , thick, draw, inner sep=2.5]{};
\node(11) at (11,0) [circle, blue   , thick, draw, inner sep=2.5]{};

\draw[thick] (6)  -- (7);
\draw[thick] (7)  -- (8);
\draw[thick,dashed] (8) -- (9);
\draw[thick] (9)  -- (10);
\draw[thick] (10) -- (11);

\end{tikzpicture}
\end{equation}where a sum is to be understood on any gray node, according to the
two possible types of nodes, hollow or filled. Computing the weights
of such paths, we shall get the expressions
\begin{align}
\beta_{n}^{k}\left(u_{0}|u_{1},\ldots,u_{n}\right) & =\alpha\left(u_{k}\right)\sum_{\eta_{1},\ldots,\eta_{k-1}=1}^{2}\prod_{i=1}^{k-1}a_{\eta_{i}}\left(\lambda_{i-1}^{P},u_{i}\right)a_{2}\left(\lambda_{k-1}^{P},u_{k}\right)\prod_{j=k+1}^{n}a_{1}\left(u_{k},u_{j}\right)\nonumber \\
 & +\delta\left(u_{k}\right)\sum_{\eta_{1},\ldots,\eta_{k-1}=1}^{2}\prod_{i=1}^{k-1}d_{\eta_{i}}\left(\lambda_{i-1}^{P},u_{i}\right)d_{2}\left(\lambda_{k-1}^{P},u_{k}\right)\prod_{j=k+1}^{n}d_{1}\left(u_{k},u_{j}\right)=0, &  & 1\leqslant k\leqslant n.\label{BAEPold}
\end{align}

The most simple Bethe Ansatz equation is that one fixing the rapidity
$u_{1}$. In this case, we have simply: \begin{equation}
\begin{tikzpicture}
\node(t) at (0,0)[left]{$\beta_n^1\left(u_0|u_1,\ldots,u_n\right)\Psi_n\left(u_1,\ldots,u_n|u_{j}^{\times}\right)=\hspace{0.2cm}$};
\node(0) at (0,0) [rectangle,red, thick, draw, inner sep=2.5]{\textsc{a}};
\node(1) at (1,0) [circle,blue, fill=blue, thick, draw, inner sep=2.5]{};
\node(2) at (2,0) [circle,blue, thick, draw, inner sep=2.5]{};
\node(3) at (3,0) [circle,blue, thick, draw, inner sep=2.5]{};

\draw[thick] (0)--(1);
\draw[thick] (1)--(2);
\draw[thick,dashed] (2)--(3);

\node(P)  at (3.45,0){$+$};
\node(4) at (4,0) [rectangle,red, thick, draw, inner sep=2.5]{\textsc{d}};
\node(5) at (5,0) [circle,blue   ,fill=blue, thick, draw, inner sep=2.5]{};
\node(6) at (6,0) [circle,blue   , thick, draw, inner sep=2.5]{};
\node(7) at (7,0) [circle,blue   , thick, draw, inner sep=2.5]{};

\draw[thick] (4)--(5);
\draw[thick] (5)--(6);
\draw[thick,dashed] (6)--(7);

\end{tikzpicture}
\end{equation} that is, 
\begin{equation}
\beta_{n}^{1}\left(u_{0}|u_{1},\ldots,u_{n}\right)=\alpha\left(u_{1}\right)a_{2}\left(u_{0},u_{1}\right)\prod_{i=2}^{n}a_{1}\left(u_{1},u_{i}\right)+\delta\left(u_{1}\right)d_{2}\left(u_{0},u_{1}\right)\prod_{i=2}^{n}d_{1}\left(u_{1},u_{i}\right)=0,
\end{equation}
which agrees with (\ref{BAEP}) for $k=1$. The other Bethe Ansatz
equations, however, are not yet in the same form as (\ref{BAEP}).
They can, nevertheless, be simplified and then be put into the same
form as (\ref{BAEP}), after we make use of the symmetry of the excited
states regarding the permutation of the rapidities. To see that, notice
that the commutation relation (\ref{BB6VP}) implies that the $n$th
excited state does not change if any pair of rapidities are permuted.
In particular, we have that $\Psi_{n}\left(u_{1},u_{2},\ldots,u_{k},\ldots,u_{n}\right)=\Psi_{n}\left(u_{k},u_{2},\ldots,u_{1},\ldots,u_{n}\right)$,
where the rapidities $u_{1}$ and $u_{k}$ have been permuted. If
we consider the respective diagrams for the action of the diagonal
operators $A(u_{0})$ and $D(u_{0})$ on this ``not well-ordered''
state $\Psi_{n}\left(u_{k},u_{2},\ldots,u_{1},\ldots,u_{n}\right)$,
then we would find that these diagrams are identical to the original
ones, except that the labels $u_{1}$ and $u_{k}$ would be exchanged.
Since the result in both cases must be the same, we conclude that
the Bethe Ansatz equation fixing $u_{k}$ can also be written in the
same form as the Bethe Ansatz equation fixing $u_{1}$, provided that
$u_{k}$ and $u_{1}$ are permuted. 

This leads us to a simpler rule determining the Bethe Ansatz equation
for the rapidity $u_{k}$ $(2\leqslant k\leqslant n)$:
\begin{itemize}
\item The Bethe Ansatz equation fixing the rapidity $u_{k}$ $(2\leqslant k\leqslant n)$
can be obtained by the same paths that determine the Bethe Ansatz
equation fixing $u_{1}$, provided that the labels $u_{1}$ and $u_{k}$
are permuted.
\end{itemize}
We conclude therefore that the final form of the Bethe Ansatz equations
is, 
\begin{equation}
\beta_{n}^{k}\left(u_{0}|u_{1},\ldots,u_{n}\right)=\alpha\left(u_{k}\right)a_{2}\left(u_{0},u_{k}\right)\prod_{i=1,i\neq k}^{n}a_{1}\left(u_{k},u_{i}\right)+\delta\left(u_{k}\right)d_{2}\left(u_{0},u_{k}\right)\prod_{i=1,i\neq k}^{n}d_{1}\left(u_{k},u_{i}\right)=0,\quad1\leqslant k\leqslant n,\label{BAEPnew}
\end{equation}
 which is indeed equal to (\ref{BAEP}) for any value of $k$.

We highlight that the proofs of (\ref{EigenPnew}) and (\ref{BAEPnew})
follow from combinatorial arguments only \textendash{} there is no
need of using mathematical induction here. It also not necessary to
analyze the first cases (e.g., $n=1$ or $n=2$) first: the results
follow once and for all for general $n$. 

Finally, we remark that the equality between (\ref{BAEPold}) and
(\ref{BAEPnew}) also provides several intricate identities among
the elements of the $R$-matrix, namely, the following ones:
\begin{equation}
\sum_{\eta_{1},\ldots,\eta_{k-1}=1}^{2}\prod_{i=1}^{k-1}a_{\eta_{i}}\left(\lambda_{i-1}^{P},u_{i}\right)a_{2}\left(\lambda_{k-1}^{P},u_{k}\right)=a_{2}\left(u_{0},u_{k}\right)\prod_{i=1}^{k-1}a_{1}\left(u_{k},u_{i}\right),\qquad1\leqslant k\leqslant n,
\end{equation}
 and 
\begin{equation}
\sum_{\eta_{1},\ldots,\eta_{k-1}=1}^{2}\prod_{i=1}^{k-1}d_{\eta_{i}}\left(\lambda_{i-1}^{P},u_{i}\right)d_{2}\left(\lambda_{k-1}^{P},u_{k}\right)=d_{2}\left(u_{0},u_{k}\right)\prod_{i=1}^{k-1}d_{1}\left(u_{k},u_{i}\right),\qquad1\leqslant k\leqslant n.
\end{equation}
 The identities (\ref{ID6P1}) and (\ref{ID6P2}) are just special
cases of the identities above.

\section{The six-vertex model with non-periodic boundary conditions\label{Section-Boundary}}

\subsection{The monodromy and transfer matrices of the six-vertex model with
non-periodic boundary conditions \label{Section-SM-B}}

The integrability of systems described by non-periodic boundary conditions
is studied through the \emph{boundary algebraic Bethe Ansatz}. In
this case, we consider a chain with $N$ lines and $L$ columns as
before, however, we do not identify the column $L+1$ and the row
$N+1$ with the first ones, respectively. We assume, on the contrary,
that the vertices at the boundaries are different from the remaining
ones. It is convenient, therefore, to consider these boundary vertices
are disposed at the columns $j=0$ and $j=L+1$ and at the rows $i=0$
and $i=N+1$ of the lattice. These boundary vertices can also have
your own configurations, so that the respective Boltzmann weights
can be represented by two \emph{boundary matrices} \textendash{} also
known as \emph{reflection matrices} or simply $K$-matrices. Thus,
there are two $K$-matrices for each row of the lattice: we write
$K^{+}(u)$ for the left reflection matrix and $K^{-}(u)$ for the
right one. Since there is nothing beyond the boundaries, these $K$-matrices
should act only in a Hilbert space that is isomorphic to $\mathbb{C}^{2}$
\textendash{} this should be contrasted with the $R$-matrix, which
is defined on a Hilbert space isomorphic to $\mathbb{C}^{2}\otimes\mathbb{C}^{2}$.

\begin{figure}[H]
\begin{center}
\begin{tikzpicture}[ultra thick,line cap=round]

\draw [domain=90:270] plot ({ 0+cos(\x)}, {1+sin(\x)});
\draw [domain=-90:90] plot ({11  +cos(\x)}, {1+sin(\x)}); 

\draw         (0,0) -- ( 4,0);
\draw[dashed] (5,0) -- ( 6,0);
\draw         (7,0) -- (11,0);
\draw         (0,2) -- ( 4,2);
\draw[dashed] (5,2) -- ( 6,2);
\draw         (7,2) -- (11,2);

\draw (1, -0.5) -- (1 ,0.5);
\draw (3, -0.5) -- (3 ,0.5);
\draw (8, -0.5) -- (8 ,0.5);
\draw (10,-0.5) -- (10,0.5);
\draw ( 1, 2.5) -- (1 ,1.5);
\draw ( 3, 2.5) -- (3 ,1.5);
\draw ( 8, 2.5) -- (8 ,1.5);
\draw (10, 2.5) -- (10,1.5);

\filldraw[red] (-1,1) circle  (3pt);
\filldraw[red] (12,1) circle  (3pt);

\filldraw[blue] (1 ,0) circle  (3pt);
\filldraw[blue] (3 ,0) circle  (3pt);
\filldraw[blue] (8 ,0) circle  (3pt);
\filldraw[blue] (10,0) circle (3pt);
6

\filldraw[blue] (1 ,2) circle  (3pt);
\filldraw[blue] (3 ,2) circle  (3pt);
\filldraw[blue] (8 ,2) circle  (3pt);
\filldraw[blue] (10,2) circle  (3pt);

\node[left]  at (-1.1 ,1){$K_a^{+}(u)$};
\node[right] at (12.1 ,1){$K_a^{-}(u)$};

\node[below] at (1 ,-0.5){$R_{a1}(u)$};
\node[below] at (3 ,-0.5){$R_{a2}(u)$};
\node[below] at (8 ,-0.5){$R_{a(L-1)}(u)$};
\node[below] at (10,-0.5){$R_{aL}(u)$};

\node[above] at (1 ,2.5){$R^{-1}_{a1}(-u)$};
\node[above] at (3 ,2.5){$R^{-1}_{a2}(-u)$};
\node[above] at (8 ,2.5){$R^{-1}_{a(L-1)}(-u)$};
\node[above] at (10 ,2.5){$R^{-1}_{aL}(-u)$};

\end{tikzpicture} 
\end{center}

\caption{A graphical representation of the double monodromy. This is defined
by the product of all $R$-matrices in a given row ($i=a$) of the
lattice, taken in the usual order, times the $K$-matrix at column
$j=L+1$, times the product of the inverse of these $R$-matrices
with the opposite spectral parameter and in the reversed order. All
these products follow the anti clock-wise direction of the diagram
above. If we further multiply this quantity by the other $K$-matrix
in column $j=0$ of the lattice and then we take the trace over the
vector space $H_{a}$, we shall obtain the transfer matrix.}

\label{Fig-DoubleB}
\end{figure}
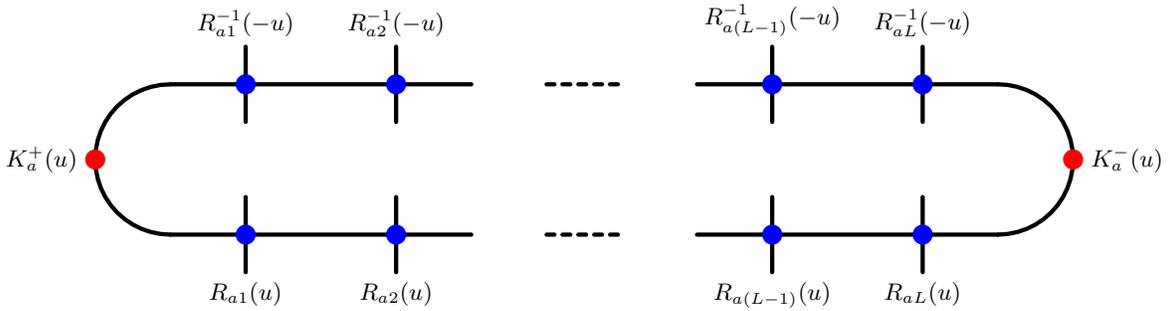

In the non-periodic case, we cannot define a single monodromy matrix
for a given row of the lattice because we need to take account for
contributions at the boundaries. One way of working this around is
to consider a \emph{double monodromy}, which was introduced by Sklyanin
\cite{Sklyanin1988}. This double monodromy matrix is constructed
as follows: we first multiply all $R$-matrices in a given row ($i=a$)
of the lattice \textendash{} that is, we consider the usual periodic
monodromy, $M_{a}(u)$ \textendash , then we multiply it by $K$-matrix
at the right boundary and, finally, we return back by multiplying
it by the inverse of all those $R$-matrices in the reversed order
and with the opposite of the spectral parameter \textendash{} that
is, we multiply by the inverse of the periodic monodromy, $M_{a}^{-1}(-u)$,
which we might call the \emph{reflected monodromy}. That is, the double
monodromy $U_{a}(u)$ is therefore defined as follows: 
\begin{equation}
U_{a}(u)=M_{a}(u)K^{-}(u)M_{a}^{-1}(-u),\label{Ua}
\end{equation}
where, 
\begin{equation}
M_{a}(u)=R_{a1}(u)\cdots R_{aL}(u),\qquad\text{and}\qquad M_{a}^{-1}(-u)=R_{aL}^{-1}(-u)\cdots R_{a1}^{-1}(-u)=\frac{R_{a1}(u)\cdots R_{aL}(u)}{r_{1}^{L}(u)r_{1}^{L}(-u)}.\label{MM}
\end{equation}
The last identity is obtained from the relation $R_{ij}^{-1}(-u)=R_{ij}(u)\left/r_{1}(u)r_{1}(-u)\right.$,
which holds for the symmetric six-vertex $R$-matrix (\ref{R0}).
See Fig. \ref{Fig-DoubleB} for a graphical representation of the
double monodromy.

After the double monodromy matrix (\ref{Ua}) is defined, the transfer
matrix for the non-periodic case can be constructed. This is given
by the trace (on the vector space $H_{a}$) of the boundary matrix
$K^{+}(u)$ multiplied by the double monodromy $U_{a}(u)$: 
\begin{equation}
T(u)=\mathrm{tr}_{a}\left[K^{+}(u)U_{a}(u)\right].
\end{equation}

The partition function of the whole lattice can then be found by repeating
the above procedure for all rows of the lattice, as we take into account
the boundary matrices at the beginning and the ending of each column.
When all the rows in the bulk of the lattice are equivalent, the eigenvalue
for the transfer matrix of any row will be the same. Thus, as in the
periodic case, the eigenvalues of the partition function will depend
only on the eigenvalues of the transfer matrix, so that we can concern
ourselves only with the diagonalization of the transfer matrix.

\subsection{The algebraic Bethe Ansatz \label{Section-ABA-B}}

The algebraic Bethe Ansatz can be generalized to cover the case where
non-periodic boundary conditions take place. In this case, we call
it the boundary algebraic Bethe Ansatz. The steps necessary to execute
the algebraic Bethe Ansatz, that is, to find the eigenvalues and eigenstates
of the transfer matrix, are similar to that of the periodic case,
although it is more complex due to the existence of the boundary matrices.
We list below the main steps:
\begin{enumerate}
\item \emph{The solutions of the boundary Yang-Baxter equations providing
the reflection $K$-matri}ces. Assuming that there is known an $R$-matrix,
solution of the Yang-Baxter equation, we look for the reflection matrices
$K^{\pm}(u)$. These matrices are the solutions of the boundary Yang-Baxter
equations that ensure the integrability at the boundaries; 
\item \emph{The Lax representation for the double monodromy and transfer
matrices. }As in the periodic case, we implement a representation
in which the double monodromy explicitly exhibit the annihilator and
creator operators used in the construction of the excited states.
The transfer matrix also become given by the diagonal operators of
the monodromy matrix, times the elements of the left boundary matrix; 
\item \emph{The reference state.} It has the same meaning as in the periodic
case;
\item \emph{The construction of the excited states}. Built in the same way
as in the periodic case, that is, through the action of the creator
operators on the reference state;
\item \emph{The derivation of the commutation relations}. They are obtained
from the fundamental relation of the boundary algebraic Bethe Ansatz;
\item \emph{The computation of the eigenvalues of the transfer matrix}.
This is again the hardest step. To compute the eigenvalues we need
to use the commutation relations repeatedly;
\item \emph{The solution of the Bethe Ansatz equations}. The consistency
conditions for the transfer matrix satisfy the eigenvalue equation.
\end{enumerate}
We shall explain these steps in details in the following.

\subsubsection{The solutions of the boundary Yang-Baxter equations, the boundary
$K$-matrices}

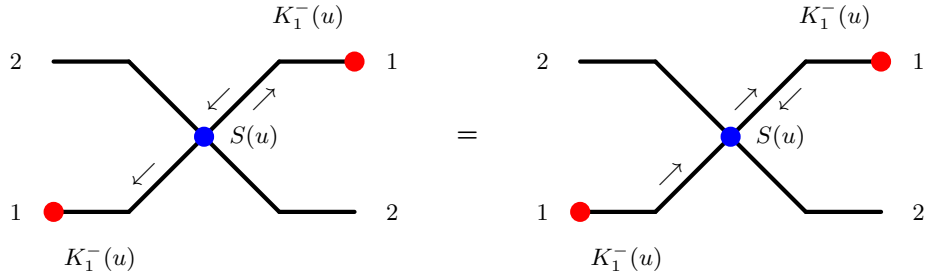
\begin{figure}[H]
\begin{center}
\begin{tikzpicture}[ultra thick,line cap=round]
\draw (-2,-1) -- (-1, -1);
\draw (-2, 1) -- (-1,  1);
\draw (-1,-1) -- (1, 1);
\draw (-1, 1) -- (1,-1);
\draw ( 1, 1) -- (2, 1);
\draw ( 1,-1) -- (2,-1);

\node at (3.5,0){$\bm{=}$};

\draw ( 5, 1) -- (6, 1);
\draw ( 5,-1) -- (6,-1);
\draw ( 6,-1) -- (8, 1);
\draw ( 6, 1) -- (8,-1);
\draw (8 ,-1) -- (9, -1);
\draw (8 , 1) -- (9,  1);

\node[right] at ( 0.5, 0.5){$\nearrow$};
\node[left]  at ( 0.5, 0.5){$\swarrow$};
\node[left]  at (-0.5,-0.5){$\swarrow$};

\node[right]  at (7.5, 0.5){$\swarrow$};
\node[left]   at (7.5, 0.5){$\nearrow$};
\node[left]   at (6.5,-0.5){$\nearrow$};

\filldraw[blue] ( 0, 0) circle  (3pt);
\filldraw[blue] ( 7, 0) circle  (3pt);
\filldraw[ red] ( 2, 1) circle  (3pt);]
\filldraw[ red] (-2,-1) circle  (3pt);
\filldraw[ red] ( 5,-1) circle  (3pt);
\filldraw[ red] ( 9, 1) circle  (3pt);

\node at (-2.5, 1){$2$};
\node at (-2.5, -1){$1$};
\node at (2.5, 1){$1$};
\node at (2.5, -1){$2$};
\node at (4.5, 1){$2$};
\node at (4.5, -1){$1$};
\node at (9.5, 1){$1$};
\node at (9.5,-1){$2$};

\node[right]  at (0.2,   0){$S(u)$};
\node[right]  at (7.2,   0){$S(u)$};
\node[above left] at (   2, 1.3){$K_{1}^{-}(u)$};
\node[below right] at (  -2,-1.3){$K_{1}^{-}(u)$};
\node[above left]  at (  9, 1.3){$K_{1}^{-}(u)$};
\node[below right] at (  5,-1.3){$K_{1}^{-}(u)$};

\end{tikzpicture} 
\end{center}

\caption{A graphical representation of the boundary Yang-Baxter Equation (\ref{BYBS1}).
The $K$-matrices behave as a reflection in the horizontal direction.
The arrows indicate the order in which the left and right diagrams
should be read, starting from the center in the first diagram and
from the upper left on the second. }

\label{Fig-BYB}
\end{figure}

The integrability of the system at the boundaries is guaranteed the
\emph{boundary Yang-Baxter equations} \textendash{} also known as
the \emph{reflection equations} \textendash{} which are \cite{Sklyanin1988},
\begin{equation}
R\left(u-v\right)K_{1}^{-}(u)R(u+v)K_{2}^{-}(v)=K_{2}^{-}(v)R(u+v)K_{1}^{-}(u)R\left(u-v\right),\label{BYBE1}
\end{equation}
and
\begin{equation}
R\left(v-u\right)\left[K_{1}^{+}(u)\right]^{t_{1}}R\left(-u-v-2\rho\right)\left[K_{2}^{+}(v)\right]^{t_{2}}=\left[K_{2}^{+}(v)\right]^{t_{2}}R\left(-u-v-2\rho\right)\left[K_{1}^{+}(u)\right]^{t_{1}}R\left(v-u\right).\label{BYBE2}
\end{equation}
In terms of the $S$ matrix (\ref{S}), the reflection equations above
can be written, respectively, as, 
\begin{equation}
S\left(u-v\right)K_{1}^{-}(u)S(u+v)K_{1}^{-}(v)=K_{1}^{-}(v)S(u+v)K_{1}^{-}(u)S\left(u-v\right),\label{BYBS1}
\end{equation}
 and 
\begin{equation}
S\left(v-u\right)\left[K_{1}^{+}(u)\right]^{t_{1}}S\left(-u-v-2\rho\right)\left[K_{1}^{+}(v)\right]^{t_{2}}=\left[K_{1}^{+}(v)\right]^{t_{2}}S\left(-u-v-2\rho\right)\left[K_{1}^{+}(u)\right]^{t_{1}}S\left(v-u\right),\label{BYBS2}
\end{equation}
 which are more useful to be explained graphically \textendash{} see
Fig. \ref{Fig-BYB}. 

The reflection equations for $K^{-}$ ensure the integrability at
the right side of the lattice, while the reflection equations for
$K^{+}$ guarantee the integrability at the left. The integrability
in the bulk is still provided by the periodic Yang-Baxter equation
(\ref{YBE}). The reflection equations are defined in $\mathrm{End}\left(V\otimes V\right)$,
where $V$ is isomorphic to $\mathbb{C}^{2}$. The reflection $K$-matrices,
by they turn, are defined in $\mathrm{End}\left(V\right)$ so that
$K_{1}^{\pm}=K^{\pm}\otimes I$ and $K_{2}^{\pm}=I\otimes K^{\pm}$,
with $I$ denoting the identity matrix belonging to $\mathrm{End}\left(V\right)$.
Besides, $t_{1}$ and $t_{2}$ mean the partial transposition in the
first and second vector spaces, respectively, and $\rho$ is the so-called
\emph{crossing parameter} \textendash{} a parameter specific to the
model that provides, for instance, the isomorphism $K^{+}(u)=\left[K^{-}(-u-\rho)\right]^{t}$
between the solutions of (\ref{BYBE1}) and (\ref{BYBE2}) \textendash{}
see \cite{Sklyanin1988} for more details.

For the six-vertex model, the $K$-matrices are two-by-two matrices.
We are actually only interested in the diagonal reflection $K$-matrices
of the six-vertex-model, which are of the form,
\begin{equation}
K^{\pm}(u)=\begin{pmatrix}k_{1,1}^{\pm}(u) & 0\\
0 & k_{2,2}^{\pm}(u)
\end{pmatrix}.
\end{equation}

The methods of solving (\ref{BYBE1}) and (\ref{BYBE2}) are somewhat
similar to that employed in the periodic case, so that we shall despise
it. Notwithstanding, there are several solutions for the six-vertex
model with diagonal boundary conditions that rely for different integrable
models \textendash{} as in the periodic case again. The most important
ones are the solution associated with the \noun{xxz} Heisenberg chain
\cite{Vega1993},
\begin{align}
k_{1,1}^{-}(u) & =\sinh(u+\zeta^{-}), & k_{2,2}^{-}(u) & =-\sinh(u-\zeta^{-}),\nonumber \\
k_{1,1}^{+}(u) & =\sinh(-u-\xi+\zeta^{+}), & k_{2,2}^{+}(u) & =-\sinh(-u-\xi-\zeta^{+}),
\end{align}
and that for the \noun{xxx} Heisenberg chain, 
\begin{align}
k_{1,1}^{-}(u) & =u+\zeta^{-}, & k_{2,2}^{-}(u) & =-u+\zeta^{-},\nonumber \\
k_{1,1}^{+}(u) & =-u-\xi+\zeta^{+}, & k_{2,2}^{+}(u) & =u+\xi+\zeta^{+},
\end{align}
where $\xi^{\pm}$ and $\zeta^{\pm}$ are arbitrary parameters. 

\subsubsection{The Lax representation of the monodromy and transfer matrices}

We continue by setting up the Lax representation, that is, a representation
in which the monodromy matrix is written as a two-by-two operator-valued
matrix:
\begin{equation}
U_{a}(u)=M_{a}(u)K^{-}(u)M_{a}^{-1}(-u)=\begin{pmatrix}A(u) & B(u)\\
C(u) & D(u)
\end{pmatrix}.\label{U}
\end{equation}
 An explicit formula for the elements of $U_{a}(u)$ can be found
through (\ref{MM}) and (\ref{Mij}). 

The transfer matrix, by its turn, is given by the trace of the monodromy
matrix (\ref{U}), previously multiplied by $K^{+}(u)$, that is,
\begin{equation}
T(u)=\mathrm{tr}_{a}\left[K^{+}(u)U_{a}(u)\right]=k_{1,1}^{+}(u)A(u)+k_{2,2}^{+}(u)D(u),\label{transferB}
\end{equation}
where we considered only the case of diagonal $K$-matrices.

\subsubsection{The reference state}

Here as well, to diagonalize the transfer matrix (\ref{transferB})
through the boundary algebraic Bethe Ansatz, we shall need a reference
state. Fortunately, the same reference state found in the periodic
case holds in this case as well, namely,
\begin{equation}
\Psi_{0}=\begin{pmatrix}1\\
0
\end{pmatrix}_{1}\otimes\cdots\otimes\begin{pmatrix}1\\
0
\end{pmatrix}_{L}.\label{PsiB}
\end{equation}
In fact, it can be verified directly from the application of (\ref{MM})
that the action of the boundary monodromy elements on $\Psi_{0}$
reads, 
\begin{align}
A(u)\Psi_{0} & =\alpha(u)\Psi_{0}, & B(u)\Psi_{0} & \neq z\Psi_{0}, & C(u)\Psi_{0} & =0\Psi_{0}, & D(u)\Psi_{0} & =\delta(u)\Psi_{0},
\end{align}
for any $z\in\mathbb{C}$. In the non-periodic case, however, $\alpha(u)$
and $\delta(u)$ are given by more complicated expressions:
\begin{equation}
\alpha(u)=k_{1,1}^{-}(u)\frac{r_{1}^{2L}(u)}{r_{1}^{L}(u)r_{1}^{L}(-u)},\qquad\delta(u)=k_{2,2}^{-}(u)\frac{r_{2}^{2L}(u)}{r_{1}^{L}(u)r_{1}^{L}(-u)}+k_{1,1}^{-}(u)h_{L-1}\left(r_{1}^{2}(u),r_{2}^{2}(u)\right)\frac{r_{3}^{2}(u)}{r_{1}^{L}(u)r_{1}^{L}(-u)},\label{alphadelta1}
\end{equation}
where $h_{L}(u,v)$ denotes the \emph{complete homogeneous symmetric
polynomials} of degree $L$ in two variables: 
\begin{equation}
h_{L}(u,v)=\sum_{k=0}^{L}u^{k}v^{L-k}=\frac{u^{L+1}-v^{L+1}}{u-v}.\label{HL}
\end{equation}
The action of the diagonal operators $A$ and $D$ on the reference
state $\Psi_{0}$ can also be found in another way. In fact, we could
had started with the fundamental relation of the periodic algebraic
Bethe Ansatz (\ref{FR}) and evaluated it at $v=-u$, in order get
the equation,
\begin{equation}
M_{b}^{-1}(-u)R_{ab}(2u)M_{a}(u)=M_{a}(u)R_{ab}(2u)M_{b}^{-1}(-u).
\end{equation}
From this, the commutation relations between the operators $M(u)$
and $M^{-1}(-u)$ could be obtained and, after the action of the monodromy
elements on the reference state $\Psi_{0}$ is computed, we would
find the expressions,
\begin{equation}
\alpha(u)=k_{1,1}^{-}(u)\frac{r_{1}^{2L}(u)}{r_{1}^{L}(u)r_{1}^{L}(-u)},\qquad\delta(u)=f(u)\alpha(u)+[k_{2,2}^{-}(u)-f(u)k_{1,1}^{-}(u)]\frac{r_{2}^{2L}(u)}{r_{1}^{L}(u)r_{1}^{L}(-u)},\label{alphadelta2}
\end{equation}
where 
\begin{equation}
f(u)=\frac{r_{3}(2u)}{r_{1}(2u)}=\frac{r_{3}^{2}(u)}{r_{1}^{2}(u)-r_{2}^{2}(u)}.\label{F}
\end{equation}
The equivalence between this result and that given by (\ref{alphadelta1})
is found after we use (\ref{HL}).

Therefore, the action of the transfer matrix on the reference state
reads,
\begin{equation}
T(u)\Psi_{0}=\tau_{0}(u)\Psi_{0},\qquad\tau_{0}(u)=k_{1,1}^{+}(u)\alpha(u)+k_{2,2}^{+}(u)\delta(u)
\end{equation}
where $\alpha(u)$ and $\delta(u)$ can be either written as in (\ref{alphadelta1})
or in (\ref{alphadelta2}). This proves that $\Psi_{0}$ given at
(\ref{PsiB}) is indeed an eigenstate of the transfer matrix (\ref{transferB}).

\subsubsection{The construction of the excited states}

The excited states can also be constructed in the same way as in the
periodic case, namely, through the repeated action of the creator
operator $B$ on the reference state $\Psi_{0}$:
\begin{equation}
\Psi_{n}\left(u_{1},\ldots,u_{n}\right)=B(u_{1})\cdots B(u_{n})\Psi_{0}.
\end{equation}
Hence, the action of the transfer matrix on the $n$th excited state
reads: 
\begin{equation}
T(u)\Psi_{n}\left(u_{1},\ldots,u_{n}\right)=k_{1,1}^{+}(u)A(u)B(u_{1})\cdots B(u_{n})\Psi_{0}+k_{2,2}^{+}(u)D(u)B(u_{1})\cdots B(u_{n})\Psi_{0}.\label{TPsiN-B}
\end{equation}

\subsubsection{\emph{The commutation relations}}

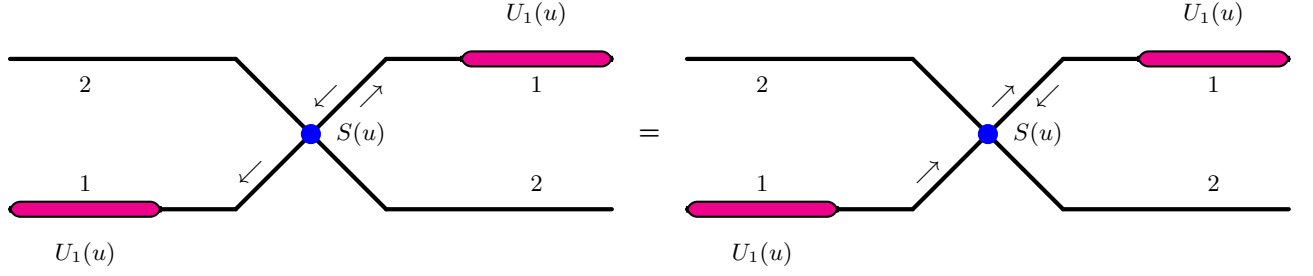
\begin{figure}[H]
\begin{center}
\begin{tikzpicture}[ultra thick,line cap=round]
\draw (-4,-1) -- (-1, -1);
\draw (-4, 1) -- (-1,  1);
\draw (-1,-1) -- (1, 1);
\draw (-1, 1) -- (1,-1);
\draw ( 1, 1) -- (4, 1);
\draw ( 1,-1) -- (4,-1);

\node at (4.5,0){$\bm{=}$};

\draw ( 5, 1) -- (8, 1);
\draw ( 5,-1) -- (8,-1);
\draw ( 8,-1) -- (10, 1);
\draw ( 8, 1) -- (10,-1);
\draw (10 ,-1) -- (13, -1);
\draw (10 , 1) -- (13,  1);

\node[right] at ( 0.5, 0.5){$\nearrow$};
\node[left]  at ( 0.5, 0.5){$\swarrow$};
\node[left]  at (-0.5,-0.5){$\swarrow$};

\node[right]  at (9.5, 0.5){$\swarrow$};
\node[left]   at (9.5, 0.5){$\nearrow$};
\node[left]   at (8.5,-0.5){$\nearrow$};

\filldraw[blue] ( 0, 0) circle  (3pt);
\filldraw[blue] ( 9, 0) circle  (3pt);
\draw[rounded corners,fill=magenta,thick] (-4,-1.1) rectangle (-2,-0.9);
\draw[rounded corners,fill=magenta,thick] ( 2, 0.9) rectangle ( 4, 1.1);
\draw[rounded corners,fill=magenta,thick] ( 5, -1.1) rectangle ( 7,-0.9);
\draw[rounded corners,fill=magenta,thick] (11,  0.9) rectangle (13,1.1);

\node[below] at ( -3, 0.9){$2$};
\node[above] at ( -3,-0.9){$1$};
\node[below] at (  3, 0.9){$1$};
\node[above] at (  3,-0.9){$2$};
\node[below] at (  6, 0.9){$2$};
\node[above] at (  6,-0.9){$1$};
\node[below] at ( 12, 0.9){$1$};
\node[above] at ( 12,-0.9){$2$};

\node[right]  at (0.2,   0){$S(u)$};
\node[right] at ( 9.2, 0){$S(u)$};
\node[above] at ( 3, 1.3){$U_{1}(u)$};
\node[below] at (-3, -1.3){$U_{1}(u)$};
\node[below] at ( 6, -1.3){$U_{1}(u)$};
\node[above] at ( 12, 1.3){$U_{1}(u)$};

\end{tikzpicture} 
\end{center}

\caption{A graphical representation of the fundamental relation of the boundary
algebraic Bethe Ansatz (\ref{BFRS}). The double monodromy is represented
by rounded rectangles. The arrows indicate the order in which the
left and right diagrams should be read, starting from the center in
the first diagram and from the upper left on the second. The fundamental
relation ensures the commutativity of the transfer matrix in the non-periodic
case, and it also provides the commutation relations between the elements
of the double monodromy, which are used in the computation of the
transfer matrix eigenvalues.}

\label{Fig-FRB}
\end{figure}

In order to compute (\ref{TPsiN-B}), the commutation relations between
the diagonal operators $A(u)$ and $D(u)$ with the creator operators
$B(u_{k})$ $(1\leqslant k\leqslant n)$ are necessary. These commutation
relations are provided by the \emph{fundamental relation of the boundary
algebraic Bethe Ansatz}: 
\begin{equation}
R\left(u-v\right)U_{a}(u)R(u+v)U_{b}(v)=U_{b}(v)R(u+v)U_{a}(u)R\left(u-v\right),\label{BFR}
\end{equation}
 or, in terms of the $S$ matrix (\ref{S}), 
\begin{equation}
S\left(u-v\right)U_{a}(u)S(u+v)U_{a}(v)=U_{a}(v)S(u+v)U_{a}(u)S\left(u-v\right).\label{BFRS}
\end{equation}
 A graphical interpretation of the fundamental relation of the boundary
algebraic Bethe Ansatz is given in Fig. \ref{Fig-FRB}.

For the six-vertex model, the needed commutation relations are the
following: 
\begin{align}
A(u)B(v) & =a_{1}(u,v)B(v)A(u)+a_{2}(u,v)B(u)A(v)+a_{3}(u,v)B(u)D(v),\label{ABold}\\
D(u)B(v) & =d_{1}(u,v)B(v)D(u)+d_{2}(u,v)B(u)D(v)+d_{3}(u,v)B(u)A(v)+d_{4}(u,v)B(v)A(u),\label{DBold}\\
B(u)B(v) & =B(v)B(u),\label{BBB}
\end{align}
where, 
\begin{align}
a_{1}(u,v) & =\frac{r_{1}(v-u)r_{2}(u+v)}{r_{1}(u+v)r_{2}(v-u)}, & a_{2}(u,v) & =-\frac{r_{2}(u+v)r_{3}(v-u)}{r_{1}(u+v)r_{2}(v-u)} & a_{3}(u,v) & =-\frac{r_{3}(u+v)}{r_{1}(u+v)},
\end{align}
 and 
\begin{align}
d_{1}(u,v) & =\frac{r_{1}(u-v)r_{1}(u+v)}{r_{2}(u-v)r_{2}(u+v)}-\frac{r_{1}(u-v)r_{3}^{2}(u+v)}{r_{1}(u+v)r_{2}(u-v)r_{2}(u+v)},\\
d_{2}(u,v) & =\frac{r_{3}(u-v)r_{3}^{2}(u+v)}{r_{1}(u+v)r_{2}(u-v)r_{2}(u+v)}-\frac{r_{1}(u+v)r_{3}(u-v)}{r_{2}(u-v)r_{2}(u+v)},\\
d_{3}(u,v) & =\frac{r_{3}(u+v)r_{1}^{2}(u-v)}{r_{1}(u+v)r_{2}^{2}(u-v)}+\frac{r_{3}(u-v)r_{3}(v-u)r_{3}(u+v)}{r_{1}(u+v)r_{2}(u-v)r_{2}(v-u)},\\
d_{4}(u,v) & =-\frac{r_{1}(u-v)r_{3}(u-v)r_{3}(u+v)}{r_{1}(u+v)r_{2}^{2}(u-v)}-\frac{r_{1}(v-u)r_{3}(u-v)r_{3}(u+v)}{r_{1}(u+v)r_{2}(u-v)r_{2}(v-u)}.
\end{align}

Notice that these expressions are more complicated than the periodic
ones. Moreover, the commutation relation between the operators $D(u)$
and $B(v)$ has one term with no counterpart on the commutation relation
between $A(u)$ and $B(v)$ \textendash{} to be more specific, the
term $d_{4}(u,v)B(v)A(u)$. This asymmetry leads to some complication
in the execution of the boundary algebraic Bethe Ansatz and, by this
reason, we usually introduce the following \emph{shifted diagonal
operators},
\begin{equation}
\mathcal{A}(u)=A(u),\qquad\mathcal{D}(u)=D(u)-f(u)A(u),\qquad f(u)=\frac{r_{3}(2u)}{r_{1}(2u)}=\frac{r_{3}^{2}(u)}{r_{1}^{2}(u)-r_{2}^{2}(u)}.
\end{equation}
order to vanish that undesirable term in (\ref{DBold}). Notice that
$f(u)$ is the same function as given by (\ref{F}). 

It follows thus that the commutation relations (\ref{ABold}), (\ref{DBold})
and (\ref{BBB}) can be brought into a more symmetric form, namely,
\begin{align}
\mathcal{A}(u)B(v) & =\mathsf{a}_{1}(u,v)B(v)\mathcal{A}(u)+\mathsf{a}_{2}(u,v)B(u)\mathcal{A}(v)+\mathsf{a}_{3}(u,v)B(u)\mathcal{D}(v),\label{AB2}\\
\mathcal{D}(u)B(v) & =\mathsf{d}_{1}(u,v)B(v)\mathcal{D}(u)+\mathsf{d}_{2}(u,v)B(u)\mathcal{D}(v)+\mathsf{d}_{3}(u,v)B(u)\mathcal{A}(v),\label{DB2}\\
B(u)B(v) & =B(v)B(u),\label{BB2}
\end{align}
where, now,
\begin{align}
\mathsf{a}_{1}(u,v) & =\frac{r_{1}(v-u)r_{2}(u+v)}{r_{1}(u+v)r_{2}(v-u)}, & \mathsf{a}_{2}(u,v) & =-\frac{r_{3}(u+v)}{r_{1}(u+v)}f(v)-\frac{r_{2}(u+v)r_{3}(v-u)}{r_{1}(u+v)r_{2}(v-u)}, & \mathsf{a}_{3}(u,v) & =-\frac{r_{3}(u+v)}{r_{1}(u+v)},
\end{align}
and,
\begin{align}
\mathsf{d}_{1}(u,v) & =\frac{r_{1}(u-v)r_{1}(u+v)}{r_{2}(u-v)r_{2}(u+v)}-\frac{r_{1}(u-v)r_{3}^{2}(u+v)}{r_{1}(u+v)r_{2}(u-v)r_{2}(u+v)},\\
\mathsf{d}_{2}(u,v) & =\frac{r_{3}(u+v)}{r_{1}(u+v)}f(u)+\frac{r_{3}(u-v)r_{3}^{2}(u+v)}{r_{1}(u+v)r_{2}(u-v)r_{2}(u+v)}-\frac{r_{1}(u+v)r_{3}(u-v)}{r_{2}(u-v)r_{2}(u+v)},\\
\mathsf{d}_{3}(u,v) & =\frac{r_{3}(u+v)r_{1}^{2}(u-v)}{r_{1}(u+v)r_{2}^{2}(u-v)}+\frac{r_{3}(u-v)r_{3}(v-u)r_{3}(u+v)}{r_{1}(u+v)r_{2}(u-v)r_{2}(v-u)}+\left[\frac{r_{2}(u+v)r_{3}(v-u)}{r_{1}(u+v)r_{2}(v-u)}+\frac{r_{3}(u+v)}{r_{1}(u+v)}f(v)\right]f(u)\nonumber \\
 & +\left[\frac{r_{3}(u-v)r_{3}^{2}(u+v)}{r_{1}(u+v)r_{2}(u-v)r_{2}(u+v)}-\frac{r_{1}(u+v)r_{3}(u-v)}{r_{2}(u-v)r_{2}(u+v)}\right]f(v).
\end{align}

The action of the shifted diagonal operators in the reference state
then becomes,
\begin{align}
\mathfrak{a}(u) & =\alpha(u)=k_{1,1}^{-}(u)\frac{r_{1}^{2L}(u)}{r_{1}^{L}(u)r_{1}^{L}(-u)},\\
\mathfrak{d}(u) & =\delta(u)-f(u)\alpha(u)\nonumber \\
 & =k_{2,2}^{-}(u)\frac{r_{2}^{2L}(u)}{r_{1}^{L}(u)r_{1}^{L}(-u)}+k_{1,1}^{-}(u)h_{L-1}\left(r_{1}^{2}(u),r_{2}^{2}(u)\right)\frac{r_{3}^{2}(u)}{r_{1}^{L}(u)r_{1}^{L}(-u)}-f(u)k_{1,1}^{-}(u)\frac{r_{1}^{2L}(u)}{r_{1}^{L}(u)r_{1}^{L}(-u)}\nonumber \\
 & =[k_{2,2}^{-}(u)-f(u)k_{1,1}^{-}(u)]\frac{r_{2}^{2L}(u)}{r_{1}^{L}(u)r_{1}^{L}(-u)}.
\end{align}

Finally, introducing as well the shifted $K$-matrices elements, 
\begin{equation}
\kappa_{1,1}^{+}(u)=k_{1,1}^{+}(u)+f(u)k_{2,2}^{+}(u),\qquad\kappa_{2,2}^{+}(u)=k_{2,2}^{+}(u),
\end{equation}
the transfer matrix can be rewritten as, 
\begin{equation}
T(u)=\kappa_{1,1}^{+}(u)\mathcal{A}(u)+\kappa_{2,2}^{+}(u)\mathcal{D}(u)\Psi_{0},
\end{equation}
so that we have as well,
\begin{equation}
T(u)\Psi_{0}=\tau_{0}(u)\Psi_{0},\qquad\tau_{0}(u)=\kappa_{1,1}^{+}(u)\mathfrak{a}(u)+\kappa_{2,2}^{+}(u)\mathfrak{d}(u).\label{transferB2}
\end{equation}

\subsubsection{The computation of the eigenvalues of the transfer matrix}

At this point we have all the ingredients to execute the boundary
algebraic Bethe Ansatz. The execution is, however, very laborious,
since we have to use the commutation relations (\ref{AB2}) and (\ref{DB2})
several times in order to find the action of the transfer matrix on
the excited states. In fact, we usually analyzes the first cases first
\textendash{} as in the periodic case \textendash , from which the
general expressions for the eigenvalues and the Bethe Ansatz equations
can be guessed. This, eventually, provides a formula for repeated
use of the commutation relations (\ref{AB2}) and (\ref{DB2}), which
need to be proved further, of course. In fact, we have that, 
\begin{align}
\mathcal{A}(u)\prod_{k=1}^{n}B(u_{k}) & =\prod_{k=1}^{n}\mathsf{a}_{1}(u,u_{k})B(u_{k})\mathcal{A}(u)\nonumber \\
 & +\sum_{k=1}^{n}\mathsf{a}_{2}(u,u_{k})\prod_{i=1,i\neq k}^{n}\mathsf{a}_{1}(u_{k},u_{i})B(u)\mathcal{A}(u_{k})+\sum_{k=1}^{n}\mathsf{a}_{3}(u,u_{k})\prod_{i=1,i\neq k}^{n}\mathsf{d}_{1}(u_{k},u_{i})B(u)\mathcal{D}(u_{k}),\label{ABB-B}\\
\mathcal{D}(u)\prod_{k=1}^{n}B(u_{k}) & =\prod_{k=1}^{n}\mathsf{d}_{1}(u,u_{k})B(u_{k})\mathcal{D}(u)\nonumber \\
 & +\sum_{k=1}^{n}\mathsf{d}_{2}(u,u_{k})\prod_{i=1,i\neq k}^{n}\mathsf{d}_{1}(u_{k},u_{i})B(u)\mathcal{D}(u_{k})+\sum_{k=1}^{n}\mathsf{d}_{3}(u,u_{k})\prod_{i=1,i\neq k}^{n}\mathsf{a}_{1}(u_{k},u_{i})B(u)\mathcal{A}(u_{k}).\label{DBB-B}
\end{align}
These formulas can be proved through mathematical induction over $n$,
in a similar fashion as the periodic case (although the calculations
are somewhat more cumbersome). In fact, for $n=1$ this is trivially
satisfied. Suppose then that the propositions (\ref{ABB-B}) and (\ref{DBB-B})
hold for general $n$. Then, for $n+1$ we have (after we use the
commutative property of the creator operators in order to pass $B(u_{n+1})$
for the left side),
\begin{align}
\mathcal{A}(u)\prod_{k=1}^{n+1}B(u_{k}) & =\mathsf{a}_{1}(u,u_{n+1})B(u_{n+1})\mathcal{A}(u)\prod_{k=1}^{n}B(u_{k})\nonumber \\
 & +\mathsf{a}_{2}(u,u_{n+1})B(u)\mathcal{A}(u_{n+1})\prod_{k=1}^{n}B(u_{k})+\mathsf{a}_{3}(u,u_{n+1})B(u)\mathcal{D}(u_{n+1})\prod_{k=1}^{n}B(u_{k}),\\
\mathcal{D}(u)\prod_{k=1}^{n+1}B(u_{k}) & =\mathsf{d}_{1}(u,u_{n+1})B(u_{n+1})\mathcal{D}(u)\prod_{k=1}^{n}B(u_{k})\nonumber \\
 & +\mathsf{d}_{2}(u,u_{n+1})B(u)\mathcal{D}(u_{n+1})\prod_{k=1}^{n}B(u_{k})+\mathsf{d}_{3}(u,u_{n+1})B(u)\mathcal{A}(u_{n+1})\prod_{k=1}^{n}B(u_{k}),
\end{align}
 that is, 
\begin{align}
\mathcal{A}(u)\prod_{k=1}^{n+1}B(u_{k}) & =\mathsf{a}_{1}(u,u_{n+1})B(u_{n+1})\prod_{k=1}^{n}\mathsf{a}_{1}(u,u_{k})B(u_{k})\mathcal{A}(u)\nonumber \\
 & +\mathsf{a}_{1}(u,u_{n+1})B(u_{n+1})\sum_{k=1}^{n}\mathsf{a}_{2}(u,u_{k})\prod_{i=1,i\neq k}^{n}\mathsf{a}_{1}(u_{k},u_{i})B(u)\mathcal{A}(u_{k})\nonumber \\
 & +\mathsf{a}_{1}(u,u_{n+1})B(u_{n+1})\sum_{k=1}^{n}\mathsf{a}_{3}(u,u_{k})\prod_{i=1,i\neq k}^{n}\mathsf{d}_{1}(u_{k},u_{i})B(u)\mathcal{D}(u_{k})\nonumber \\
 & +\mathsf{a}_{2}(u,u_{n+1})B(u)\prod_{k=1}^{n}\mathsf{a}_{1}(u_{n+1},u_{k})B(u_{k})\mathcal{A}(u_{n+1})\nonumber \\
 & +\mathsf{a}_{2}(u,u_{n+1})B(u)\sum_{k=1}^{n}\mathsf{a}_{2}(u_{n+1},u_{k})\prod_{i=1,i\neq k}^{n}\mathsf{a}_{1}(u_{k},u_{i})B(u_{n+1})\mathcal{A}(u_{k})\nonumber \\
 & +\mathsf{a}_{2}(u,u_{n+1})B(u)\sum_{k=1}^{n}\mathsf{a}_{3}(u_{n+1},u_{k})\prod_{i=1,i\neq k}^{n}\mathsf{d}_{1}(u_{k},u_{i})B(u_{n+1})\mathcal{D}(u_{k})\nonumber \\
 & +\mathsf{a}_{3}(u,u_{n+1})B(u)\prod_{k=1}^{n}\mathsf{d}_{1}(u_{n+1},u_{k})B(u_{k})\mathcal{D}(u)\nonumber \\
 & +\mathsf{a}_{3}(u,u_{n+1})B(u)\sum_{k=1}^{n}\mathsf{d}_{2}(u_{n+1},u_{k})\prod_{i=1,i\neq k}^{n}\mathsf{d}_{1}(u_{k},u_{i})B(u_{n+1})\mathcal{D}(u_{k})\nonumber \\
 & +\mathsf{a}_{3}(u,u_{n+1})B(u)\sum_{k=1}^{n}\mathsf{d}_{3}(u_{n+1},u_{k})\prod_{i=1,i\neq k}^{n}\mathsf{a}_{1}(u_{k},u_{i})B(u_{n+1})\mathcal{A}(u_{k})\label{APB}
\end{align}
 and 
\begin{align}
\mathcal{D}(u)\prod_{k=1}^{n+1}B(u_{k}) & =\mathsf{d}_{1}(u,u_{n+1})B(u_{n+1})\prod_{k=1}^{n}\mathsf{d}_{1}(u,u_{k})B(u_{k})\mathcal{D}(u)\nonumber \\
 & +\mathsf{d}_{1}(u,u_{n+1})B(u_{n+1})\sum_{k=1}^{n}\mathsf{d}_{2}(u,u_{k})\prod_{i=1,i\neq k}^{n}\mathsf{d}_{1}(u_{k},u_{i})B(u)\mathcal{D}(u_{k})\nonumber \\
 & +\mathsf{d}_{1}(u,u_{n+1})B(u_{n+1})\sum_{k=1}^{n}\mathsf{d}_{3}(u,u_{k})\prod_{i=1,i\neq k}^{n}\mathsf{a}_{1}(u_{k},u_{i})B(u)\mathcal{A}(u_{k})\nonumber \\
 & +\mathsf{d}_{2}(u,u_{n+1})B(u)\prod_{k=1}^{n}\mathsf{d}_{1}(u_{n+1},u_{k})B(u_{k})\mathcal{D}(u_{n+1})\nonumber \\
 & +\mathsf{d}_{2}(u,u_{n+1})B(u)\sum_{k=1}^{n}\mathsf{d}_{2}(u_{n+1},u_{k})\prod_{i=1,i\neq k}^{n}\mathsf{d}_{1}(u_{k},u_{i})B(u_{n+1})\mathcal{D}(u_{k})\nonumber \\
 & +\mathsf{d}_{2}(u,u_{n+1})B(u)\sum_{k=1}^{n}\mathsf{d}_{3}(u_{n+1},u_{k})\prod_{i=1,i\neq k}^{n}\mathsf{a}_{1}(u_{k},u_{i})B(u_{n+1})\mathcal{A}(u_{k})\nonumber \\
 & +\mathsf{d}_{3}(u,u_{n+1})B(u)\prod_{k=1}^{n}\mathsf{a}_{1}(u_{n+1},u_{k})B(u_{k})\mathcal{A}(u)\nonumber \\
 & +\mathsf{d}_{3}(u,u_{n+1})B(u)\sum_{k=1}^{n}\mathsf{a}_{2}(u_{n+1},u_{k})\prod_{i=1,i\neq k}^{n}\mathsf{a}_{1}(u_{k},u_{i})B(u_{n+1})\mathcal{A}(u_{k})\nonumber \\
 & +\mathsf{d}_{3}(u,u_{n+1})B(u)\sum_{k=1}^{n}\mathsf{a}_{3}(u_{n+1},u_{k})\prod_{i=1,i\neq k}^{n}\mathsf{d}_{1}(u_{k},u_{i})B(u_{n+1})\mathcal{D}(u_{k}).\label{DPB}
\end{align}
Now, we should realize that these expressions can be simplified if
we make use of the following identities provided by the boundary Yang-Baxter
equation: 
\begin{align}
\mathsf{a}_{2}(u,u_{n+1})\mathsf{a}_{2}(u_{n+1},u_{k})+\mathsf{a}_{3}(u,u_{n+1})\mathsf{d}_{3}(u_{n+1},u_{k}) & =\mathsf{a}_{2}(u,u_{k})\mathsf{a}_{1}(u_{k},u_{n+1}),\\
\mathsf{a}_{2}(u,u_{n+1})\mathsf{a}_{3}(u_{n+1},u_{k})+\mathsf{a}_{3}(u,u_{n+1})\mathsf{d}_{2}(u_{n+1},u_{k}) & =\mathsf{a}_{3}(u,u_{k})\mathsf{d}_{1}(u_{k},u_{n+1}),\\
\mathsf{d}_{2}(u,u_{n+1})\mathsf{d}_{2}(u_{n+1},u_{k})+\mathsf{d}_{3}(u,u_{n+1})\mathsf{a}_{3}(u_{n+1},u_{k}) & =\mathsf{d}_{2}(u,u_{k})\mathsf{d}_{1}(u_{k},u_{n+1}),\\
\mathsf{d}_{2}(u,u_{n+1})\mathsf{d}_{3}(u_{n+1},u_{k})+\mathsf{d}_{3}(u,u_{n+1})\mathsf{a}_{2}(u_{n+1},u_{k}) & =\mathsf{d}_{3}(u,u_{k})\mathsf{a}_{1}(u_{k},u_{n+1}),
\end{align}
which hold for any $k$ from $1$ to $n$. Then, after we extend the
summations to $n+1$, we will be able to group the fifth and the ninth
terms in (\ref{APB}) and (\ref{DPB}) with the second one and also
the sixth and the eight terms with the third. With this, we shall
obtain the same expressions (\ref{ABB-B}) and (\ref{DBB-B}) but
with $n+1$ replacing $n$, so that we are done. 

Finally, the action of the transfer matrix on the $n$th excited state
follows after we multiply (\ref{ABB-B}) by $\kappa_{1,1}^{+}(u)$,
(\ref{DBB-B}) by $\kappa_{2,2}^{+}(u)$ and take the sum of these
two terms. This provides us with an expression quite analogous to
that of the periodic case, namely, 
\begin{equation}
T(u)\Psi_{n}\left(u_{1},\ldots,u_{n}\right)=\tau_{n}\left(u|u_{1},\ldots,u_{n}\right)\Psi_{n}\left(u_{1},\ldots,u_{n}\right)+\sum_{k=1}^{n}\beta_{n}^{k}\left(u|u_{1},\ldots,u_{n}\right)\Psi_{n}\left(u_{k}^{\times}\right).\label{ActionTB}
\end{equation}

Therefore, the requirement that the transfer matrix satisfies an eigenvalue
equation means that all the unwanted terms in (\ref{ActionTB}) should
vanish. This leads us to the eigenvalues, 
\begin{equation}
\tau_{n}\left(u|u_{1},\ldots,u_{n}\right)=\kappa_{1,1}^{+}\left(u\right)\mathfrak{a}(u)\prod_{k=1}^{n}\mathsf{a}_{1}\left(u,u_{k}\right)+\kappa_{2,2}^{+}\left(u\right)\mathfrak{d}(u)\prod_{k=1}^{n}\mathsf{d}_{1}\left(u,u_{k}\right),\label{eigen-B}
\end{equation}
and to the boundary Bethe Ansatz equations,
\begin{multline}
\beta_{n}^{k}\left(u|u_{1},\ldots,u_{n}\right)=\left[\mathsf{a}_{2}\left(u,u_{k}\right)\kappa_{1,1}^{+}\left(u\right)+\mathsf{d}_{3}\left(u,u_{k}\right)\kappa_{2,2}^{+}\left(u\right)\right]\mathfrak{a}\left(u_{k}\right)\prod_{i=1,i\neq k}^{n}\mathsf{a}_{1}\left(u_{k},u_{i}\right)\\
+\left[\mathsf{d}_{2}\left(u,u_{k}\right)\kappa_{2,2}^{+}\left(u\right)+\mathsf{a}_{3}\left(u,u_{k}\right)\kappa_{1,1}^{+}\left(u\right)\right]\mathfrak{d}\left(u_{k}\right)\prod_{i=1,i\neq k}^{n}\mathsf{d}_{1}\left(u_{k},u_{i}\right)=0,\qquad1\leqslant k\leqslant n.\label{BAE-B}
\end{multline}

\subsubsection{The solution of the Bethe Ansatz equations}

A completely exact solution for the spectral problem would require
an analytical solution of the Bethe Ansatz equations (\ref{BAE-B}).
Unfortunately, this is not possible up to date due to the high complexity
of these equations. Therefore, it only remains to resign ourselves
with the numerical approximations available. 

\subsection{The combinatorial approach\label{Section-Combinatorial-B}}

Now let us see how the combinatorial approach can be applied as well
for the case where non-periodic boundary conditions take place. Following
the same approach as in the periodic case, we can represent the commutation
relations (\ref{AB2}) and (\ref{DB2}) by the following combinatorial
diagrams\footnote{The diagrams regarding the non-shifted commutation relations (\ref{ABold})
and (\ref{DBold}) can also be drawn and the same analysis presented
in this section can be performed with them as well. We remark, however,
that in this case the diagrams representing the $A$ and $D$ operators
will no longer be symmetrical, which makes the analysis somewhat more
complicated. }: 
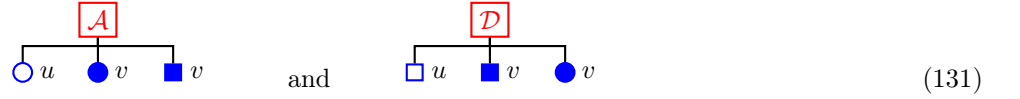
\begin{equation}
\begin{tikzpicture}[edge from parent fork down,
level 1/.style={sibling distance=1cm},
level distance=0.7cm
]
\tikzstyle{hollow circle}=[circle,blue,    thick, draw, inner sep=2.5] 
\tikzstyle{solid  circle}=[circle,blue,    thick, draw, inner sep=2.5, fill=blue]
\tikzstyle{hollow square}=[rectangle,blue, thick, draw, inner sep=3] 
\tikzstyle{solid  square}=[rectangle,blue, thick, draw, inner sep=3, fill=blue]
\tikzstyle{edge from parent}=[thick,draw]
\node(0)[hollow square,red]{$\mathcal{A}$} 
  child {node(1)[hollow circle,blue]{} }
  child {node(2)[solid  circle,blue]{} }
  child {node(3)[solid  square,blue]{} }
; 
\node[right] at (1) {$\hspace{0,1cm} u$};
\node[right] at (2) {$\hspace{0,1cm} v$};
\node[right] at (3) {$\hspace{0,1cm} v$};
\end{tikzpicture}
\hspace{1cm}
\text{and}
\hspace{1cm}
\begin{tikzpicture}[edge from parent fork down,
level 1/.style={sibling distance=1cm},
level distance=0.7cm
]
\tikzstyle{hollow circle}=[circle,,blue,    thick, draw, inner sep=2.5] 
\tikzstyle{solid  circle}=[circle,,blue,    thick, draw, inner sep=2.5, fill=blue]
\tikzstyle{hollow square}=[rectangle,blue, thick, draw, inner sep=3] 
\tikzstyle{solid  square}=[rectangle,blue, thick, draw, inner sep=3, fill=blue]
\tikzstyle{edge from parent}=[thick,draw]
\node(0)[hollow square,red]{$\mathcal{D}$} 
  child {node(1)[hollow square,blue]{} }
  child {node(2)[solid  square,blue]{} }
  child {node(3)[solid  circle,blue]{} }
; 
\node[right] at (1) {$\hspace{0,1cm} u$};
\node[right] at (2) {$\hspace{0,1cm} v$};
\node[right] at (3) {$\hspace{0,1cm} v$};
\end{tikzpicture}
\end{equation} Notice, however, that in the non-periodic case the $\mathcal{A}$
diagram has two circle nodes and one additional square node, as well
as the $\mathcal{D}$ diagram has two square nodes and one circle
node. This distinction is necessary because the third term in the
commutation relations (\ref{AB2}) and (\ref{DB2}) will lead to an
interchanging between the diagonal operators $\mathcal{A}$ and $\mathcal{D}$
when those commutation relations are applied repeatedly. In fact,
the combinatorial tree representing the action of the $\mathcal{A}$
operator on the $n$th excited state $\Psi_{n}$ will have the form:
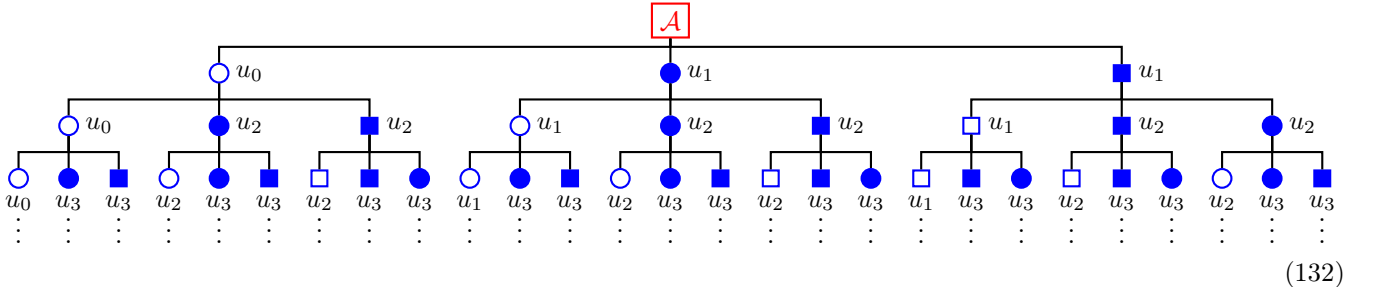
\begin{equation}
\begin{tikzpicture}[scale=1,edge from parent fork down,
level 1/.style={sibling distance=6cm},
level 2/.style={sibling distance=2cm},
level 3/.style={sibling distance=0.666cm},
level distance=0.7cm
]
\tikzstyle{hollow circle}=[circle,blue,    thick, draw, inner sep=2.5] 
\tikzstyle{solid  circle}=[circle,blue,    thick, draw, inner sep=2.5, fill=blue]
\tikzstyle{hollow square}=[rectangle,blue, thick, draw, inner sep=3] 
\tikzstyle{solid  square}=[rectangle,blue, thick, draw, inner sep=3, fill=blue]
\tikzstyle{edge from parent}=[thick,draw]

\node(0)[hollow square,red]{$\mathcal{A}$}
  child{node(1)[hollow circle]{}
    child{node(11)[hollow circle]{}
      child{node(111)[hollow circle,blue]{}}
      child{node(112)[solid  circle,blue]{}}
      child{node(113)[solid  square,blue]{}}
   }
    child{node(12)[solid circle]{}
      child{node(121)[hollow circle,blue]{}}
      child{node(122)[solid  circle,blue]{}}
      child{node(123)[solid  square,blue]{}}
   }
    child{node(13)[solid  square,blue]{}
      child{node(131)[hollow square,blue]{}}
      child{node(132)[solid square,blue]{}}
      child{node(133)[solid circle,blue]{}}
    }
  } 
  child{node(2)[solid circle,blue]{}
    child{node(21)[hollow circle,blue]{}
      child{node(211)[hollow circle,blue]{}}
      child{node(212)[solid  circle,blue]{}}
      child{node(213)[solid  square,blue]{}}
   }
    child{node(22)[solid circle,blue]{}
      child{node(221)[hollow circle,blue]{}}
      child{node(222)[solid  circle,blue]{}}
      child{node(223)[solid  square,blue]{}}
   }
    child{node(23)[solid  square,blue]{}
      child{node(231)[hollow square,blue]{}}
      child{node(232)[solid  square,blue]{}}
      child{node(233)[solid  circle,blue]{}}
    }
  } 
  child{node(3)[solid square,blue]{}
    child{node(31)[hollow square,blue]{}
      child{node(311)[hollow square,blue]{}}
      child{node(312)[solid  square,blue]{}}
      child{node(313)[solid  circle,blue]{}}
    }
    child{node(32)[solid square,blue]{}
      child{node(321)[hollow square,blue]{}}
      child{node(322)[solid  square,blue]{}}
      child{node(323)[solid  circle,blue]{}}
    }
    child{node(33)[solid circle,blue]{}
      child{node(331)[hollow circle,blue]{}}
      child{node(332)[solid circle,blue]{}}
      child{node(333)[solid  square,blue]{}}
    }
  } 
; 

\node[right=0.1cm] at  (1) {$u_0$};
\node[right=0.1cm] at  (2) {$u_1$};
\node[right=0.1cm] at  (3) {$u_1$};

\node[right=0.1cm] at  (11) {$u_0$};
\node[right=0.1cm] at  (12) {$u_2$};
\node[right=0.1cm] at  (13) {$u_2$};
\node[right=0.1cm] at  (21) {$u_1$};
\node[right=0.1cm] at  (22) {$u_2$};
\node[right=0.1cm] at  (23) {$u_2$};
\node[right=0.1cm] at  (31) {$u_1$};
\node[right=0.1cm] at  (32) {$u_2$};
\node[right=0.1cm] at  (33) {$u_2$};

\node[below=0.1cm] at  (111) {$u_0$};
\node[below=0.1cm] at  (112) {$u_3$};
\node[below=0.1cm] at  (113) {$u_3$};
\node[below=0.1cm] at  (121) {$u_2$};
\node[below=0.1cm] at  (122) {$u_3$};
\node[below=0.1cm] at  (123) {$u_3$};
\node[below=0.1cm] at  (131) {$u_2$};
\node[below=0.1cm] at  (132) {$u_3$};
\node[below=0.1cm] at  (133) {$u_3$};
\node[below=0.1cm] at  (211) {$u_1$};
\node[below=0.1cm] at  (212) {$u_3$};
\node[below=0.1cm] at  (213) {$u_3$};
\node[below=0.1cm] at  (221) {$u_2$};
\node[below=0.1cm] at  (222) {$u_3$};
\node[below=0.1cm] at  (223) {$u_3$};
\node[below=0.1cm] at  (231) {$u_2$};
\node[below=0.1cm] at  (232) {$u_3$};
\node[below=0.1cm] at  (233) {$u_3$};
\node[below=0.1cm] at  (311) {$u_1$};
\node[below=0.1cm] at  (312) {$u_3$};
\node[below=0.1cm] at  (313) {$u_3$};
\node[below=0.1cm] at  (321) {$u_2$};
\node[below=0.1cm] at  (322) {$u_3$};
\node[below=0.1cm] at  (323) {$u_3$};
\node[below=0.1cm] at  (331) {$u_2$};
\node[below=0.1cm] at  (332) {$u_3$};
\node[below=0.1cm] at  (333) {$u_3$};

\node[below=0.2cm] at  (111) {$\vdots$};
\node[below=0.2cm] at  (112) {$\vdots$};
\node[below=0.2cm] at  (113) {$\vdots$};
\node[below=0.2cm] at  (121) {$\vdots$};
\node[below=0.2cm] at  (122) {$\vdots$};
\node[below=0.2cm] at  (123) {$\vdots$};
\node[below=0.2cm] at  (131) {$\vdots$};
\node[below=0.2cm] at  (132) {$\vdots$};
\node[below=0.2cm] at  (133) {$\vdots$};
\node[below=0.2cm] at  (211) {$\vdots$};
\node[below=0.2cm] at  (212) {$\vdots$};
\node[below=0.2cm] at  (213) {$\vdots$};
\node[below=0.2cm] at  (221) {$\vdots$};
\node[below=0.2cm] at  (222) {$\vdots$};
\node[below=0.2cm] at  (223) {$\vdots$};
\node[below=0.2cm] at  (231) {$\vdots$};
\node[below=0.2cm] at  (232) {$\vdots$};
\node[below=0.2cm] at  (233) {$\vdots$};
\node[below=0.2cm] at  (311) {$\vdots$};
\node[below=0.2cm] at  (312) {$\vdots$};
\node[below=0.2cm] at  (313) {$\vdots$};
\node[below=0.2cm] at  (321) {$\vdots$};
\node[below=0.2cm] at  (322) {$\vdots$};
\node[below=0.2cm] at  (323) {$\vdots$};
\node[below=0.2cm] at  (331) {$\vdots$};
\node[below=0.2cm] at  (332) {$\vdots$};
\node[below=0.2cm] at  (333) {$\vdots$};
\end{tikzpicture}
\label{diagramAN}
\end{equation} and a similar diagram can be associated with the $\mathcal{D}$ operator
\textendash{} the only difference being that the circle nodes should
be replaced by square nodes and vice-versa.

In the diagram (\ref{diagramAN}) we adopted the convention that a
circle node means that the commutation relation between the operators
$\mathcal{A}$ and $B$ is to be used at that point further, while
a square node means that the commutation relation between $\mathcal{D}$
and $B$ is to be used at that point further. Since each term resulting
from the action of the diagonal operators on the $n$th excited state
$\Psi_{n}$ requires the use of the commutation relations $n$ times,
this will give place to a \emph{pruned ternary combinatorial trees}
of length $n$. The labels on the side of every node in the diagram
(\ref{diagramAN}) refer to the argument that the $\mathcal{A}$ operator
would have if we had used the commutation relations algebraically
(the same can be said about the $\mathcal{D}$ diagram regarding the
argument of the $\mathcal{D}$ operator, of course). These labels
can be obtained from the same rule as in the periodic case, namely,
\begin{itemize}
\item Hollow nodes inherit the label of his parent, while the label of any
filled node on the level $k$ of a given path $P$ is just $\lambda_{k}^{P}=u_{k}$
(the label of the root being $u_{0}$).
\end{itemize}
Notice that only the color (hollow or filled) of the nodes, not the
type (circle or square), determines the labels. Hence these labels
can be obtained recursively as follows: 
\begin{equation}
\lambda_{0}^{P}=u_{0},\qquad\mathrm{and}\qquad\lambda_{k}^{P}=\begin{cases}
\lambda_{k-1}^{P}, & \eta_{k}=1,\\
u_{k}, & \mathrm{otherwise},
\end{cases}\qquad\mathrm{for}1\leqslant k\leqslant n.
\end{equation}

Any path $P$ of the diagram (\ref{diagramAN}) can be uniquely specified
by the set $\eta=\left(\eta_{1},\ldots\eta_{n}\right)$, where now
the parameters $\eta_{k}$ $(1\leqslant k\leqslant n)$ can assume
the values $1$, $2$ and $3$. To any path $P$ we associate a weight
$W(P)$ and a state $\left|P\right\rangle $, so that its value $V(P)$
\textendash{} that represents the respective term obtained from the
boundary algebraic Bethe Ansatz \textendash{} can be written as $V(P)=W(P)\left|P\right\rangle $.
The weight of a given path can be defined as the product of the weights
of its nodes through the rules:
\begin{enumerate}
\item If the parent of a node in the level $k$ of a given path $P$ is
a \emph{circle}, then its weight equals $\mathsf{a}_{\eta_{k}}(\lambda_{k-1}^{P},u_{k})$;
\item If the parent of a node in the level $k$ of a given path $P$ is
a \emph{square}, then its weight equals $\mathsf{d}_{\eta_{k}}(\lambda_{k-1}^{P},u_{k})$;
\item Every leaf node in the path $P_{\mathcal{A}}$ $\left[P_{\mathcal{D}}\right]$
of the $\mathcal{A}$ $\left[\mathcal{D}\right]$ diagram also contributes
with a factor $\kappa_{1,1}^{+}(u_{0})\mathfrak{a}(\lambda_{n}^{P})\Psi_{0}$
$\left[\kappa_{2,2}^{+}(u_{0})\mathfrak{a}(\lambda_{n}^{P})\Psi_{0}\right]$
or $\kappa_{1,1}^{+}(u_{0})\mathfrak{d}(\lambda_{n}^{P})\Psi_{0}$
$\left[\kappa_{2,2}^{+}(u_{0})\mathfrak{d}(\lambda_{n}^{P})\Psi_{0}\right]$
depending on whether this leaf node is a circle or a square, respectively, 
\end{enumerate}
with the root defined as a circle for the $\mathcal{A}$ diagram and
as a square for the $\mathcal{D}$ diagram. Therefore, the weight
of any path belonging to the $\mathcal{A}$ and $\mathcal{D}$ diagrams
can be written, respectively, as\footnote{As an example, we also present the paths associated with the $\mathcal{A}$
diagram and the corresponding mathematical expressions for the second
excited state ($n=2$): 

\begin{center} %
\begin{tabular}{ll}
\begin{tikzpicture}[scale=0.7]
\node(A) at (0,0) [rectangle, thick, red, draw, inner sep=2.5]{$\mathcal{A}$}; 
\node(B) at (1,0)[circle, thick, blue, draw, inner sep=2.5]{};
\node(C) at (2,0)[circle, thick, blue, draw, inner sep=2.5]{};
\node[right] at (C) 
{\hspace{0.1cm}$ \equiv \kappa_{1,1}^{+}(u_0)a_{1}(u_{0},u_{1})a_{1}(u_{0},u_{2})\alpha\left(u_{0}\right)B(u_{1})B(u_{2})\Psi_{0},$}; 

\draw[thick] (A)--(B);
\draw[thick] (B)--(C);

\end{tikzpicture} & \begin{tikzpicture}[scale=0.7]
\node(A) at (0,0) [rectangle, thick, red, draw, inner sep=2.5]{$\mathcal{A}$}; 
\node(B) at (1,0)[circle, thick, blue, draw, inner sep=2.5]{};
\node(C) at (2,0)[circle, thick, blue, fill=blue, draw, inner sep=2.5]{};
\node[right] at (C) 
{\hspace{0.1cm}$ \equiv \kappa_{1,1}^{+}(u_0)a_{1}(u_{0},u_{1})a_{2}(u_{0},u_{2})\alpha\left(u_{2}\right)B(u_{0})B(u_{1})\Psi_{0},$};  

\draw[thick] (A)--(B);
\draw[thick] (B)--(C);

\end{tikzpicture}\tabularnewline
\begin{tikzpicture}[scale=0.7]
\tikzstyle{square}=[rectangle, blue, thick, draw, inner sep=3, fill=blue];
\node(A) at (0,0) [rectangle, thick, red, draw, inner sep=2.5]{$\mathcal{A}$}; 
\node(B) at (1,0)[circle, thick, blue, draw, inner sep=2.5]{};
\node(C) at (2,0)[square, thick, blue, draw, inner sep=3]{};
\node[right] at (C) 
{\hspace{0.1cm}$\equiv\kappa_{1,1}^{+}(u_0)a_{1}(u_{0},u_{1})d_{2}(u_{0},u_{2})\delta\left(u_{2}\right)B(u_{0})B(u_{1})\Psi_{0},$}; 
\draw[thick] (A)--(B);
\draw[thick] (B)--(C);
\end{tikzpicture} & \begin{tikzpicture}[scale=0.7]
\node(A) at (0,0) [rectangle, thick, red, draw, inner sep=2.5]{$\mathcal{A}$}; 
\node(B) at (1,0)[circle, thick, blue, fill=blue, draw, inner sep=2.5]{};
\node(C) at (2,0)[circle, thick, blue, draw, inner sep=2.5]{};
\node[right] at (C) 
{\hspace{0.1cm}$ \equiv \kappa_{1,1}^{+}(u_0)a_{2}(u_{0},u_{1})a_{1}(u_{1},u_{2})\alpha\left(u_{1}\right)B(u_{0})B(u_{2})\Psi_{0},$};
\draw[thick] (A)--(B);
\draw[thick] (B)--(C);
\end{tikzpicture}
\tabularnewline
\begin{tikzpicture}[scale=0.7]
\node(A) at (0,0) [rectangle, thick, red, draw, inner sep=2.5]{$\mathcal{A}$}; 
\node(B) at (1,0)[circle, thick, blue,fill=blue, draw, inner sep=2.5]{};
\node(C) at (2,0)[circle, thick, blue,fill=blue, draw, inner sep=2.5]{};
\node[right] at (C) 
{\hspace{0.1cm}$ \equiv \kappa_{1,1}^{+}(u_0)a_{2}(u_{0},u_{1})a_{2}(u_{1},u_{2})\alpha\left(u_{2}\right)B(u_{0})B(u_{1})\Psi_{0},$};
\draw[thick] (A)--(B);
\draw[thick] (B)--(C);

\end{tikzpicture}

 & \begin{tikzpicture}[scale=0.7]
\tikzstyle{square}=[rectangle, blue, thick, draw, inner sep=3];
\node(A) at (0,0) [rectangle, thick, red, draw, inner sep=2.5]{$\mathcal{A}$}; 
\node(B) at (1,0)[circle, thick, blue,fill=blue, draw, inner sep=2.5]{};
\node(C) at (2,0)[square, thick, blue,fill=blue, draw]{};
\node[right] at (C) 
{\hspace{0.1cm}$ \equiv \kappa_{1,1}^{+}(u_0)a_{2}(u_{0},u_{1})d_{2}(u_{1},u_{2})\delta\left(u_{2}\right)B(u_{0})B(u_{1})\Psi_{0},$};
\draw[thick] (A)--(B);
\draw[thick] (B)--(C);

\end{tikzpicture}
\tabularnewline
\begin{tikzpicture}[scale=0.7]
\tikzstyle{square}=[rectangle, blue, thick, draw, inner sep=3];
\node(A) at (0,0) [rectangle, thick, red, draw, inner sep=2.5]{$\mathcal{A}$}; 
\node(B) at (1,0)[square, thick, blue, fill=blue, draw]{};
\node(C) at (2,0)[square, thick, blue, draw]{};
\node[right] at (C) 
{\hspace{0.1cm}$ \equiv \kappa_{1,1}^{+}(u_0)d_{2}(u_{0},u_{1})d_{1}(u_{1},u_{2})\delta\left(u_{1}\right)B(u_{0})B(u_{2})\Psi_{0},$};
\draw[thick] (A)--(B);
\draw[thick] (B)--(C);
\end{tikzpicture}
 & \begin{tikzpicture}[scale=0.7]
\tikzstyle{square}=[rectangle, blue, thick, draw, inner sep=3];
\node(A) at (0,0) [rectangle, thick, red, draw, inner sep=2.5]{$\mathcal{A}$}; 
\node(B) at (1,0)[square, thick, blue,fill=blue, draw]{};
\node(C) at (2,0)[square, thick, blue,fill=blue, draw]{};
\node[right] at (C) 
{\hspace{0.1cm}$ \equiv \kappa_{1,1}^{+}(u_0)d_{2}(u_{0},u_{1})d_{2}(u_{1},u_{2})\delta\left(u_{2}\right)B(u_{0})B(u_{1})\Psi_{0},$};
\draw[thick] (A)--(B);
\draw[thick] (B)--(C);

\end{tikzpicture}
\tabularnewline
\begin{tikzpicture}[scale=0.7]
\tikzstyle{square}=[rectangle, blue, thick, draw, inner sep=3];
\node(A) at (0,0) [rectangle, thick, red, draw, inner sep=2.5]{$\mathcal{A}$}; 
\node(B) at (1,0)[square, thick, blue, fill=blue, draw]{};
\node(C) at (2,0)[circle, thick, blue, fill=blue, draw, inner sep=2.5]{};
\node[right] at (C) 
{\hspace{0.1cm}$ \equiv \kappa_{1,1}^{+}(u_0)d_{2}(u_{0},u_{1})a_{2}(u_{1},u_{2})\alpha\left(u_{2}\right)B(u_{0})B(u_{1})\Psi_{0}.$};
\draw[thick] (A)--(B);
\draw[thick] (B)--(C);
\end{tikzpicture}

 & \tabularnewline
\end{tabular}\end{center} The paths associated with the $\mathcal{D}$ diagram
are similar to the above ones, but we need to replace $\kappa_{1,1}^{+}$
by $\kappa_{2,2}^{+}$, $\alpha$ by $\delta$ (and vice-versa) and
the coefficients $a_{1}$, $a_{2}$ should be replaced by $d_{1}$,
$d_{2}$, respectively (and vice-versa).} 

\begin{align}
W\left[P_{\mathcal{A}}\left(\eta_{1},\ldots,\eta_{n}\right)\right] & =\kappa_{1,1}^{+}\left(u_{0}\right)\mathfrak{a}\left(\lambda_{n}^{P}\right)\prod_{k=1}^{n}\mathsf{b}_{\eta_{k}}^{\eta_{k-1}}\left(\lambda_{k-1}^{P},u_{k}\right),\\
W\left[P_{\mathcal{D}}\left(\eta_{1},\ldots,\eta_{n}\right)\right] & =\kappa_{2,2}^{+}\left(u_{0}\right)\mathfrak{d}\left(\lambda_{n}^{P}\right)\prod_{k=1}^{n}\mathsf{c}_{\eta_{k}}^{\eta_{k-1}}\left(\lambda_{k-1}^{P},u_{k}\right),
\end{align}
 where, 
\begin{equation}
\mathsf{b}_{\eta_{k}}^{\eta_{k-1}}=\begin{cases}
\mathsf{a}_{\eta_{k}}, & \eta_{k-1}=1,2,\\
\mathsf{d}_{\eta_{k}}, & \mathrm{otherwise},
\end{cases}\qquad\mathsf{c}_{\eta_{k}}^{\eta_{k-1}}=\begin{cases}
\mathsf{d}_{\eta_{k}}, & \eta_{k-1}=1,2,\\
\mathsf{a}_{\eta_{k}}, & \mathrm{otherwise},
\end{cases}\qquad\text{with},\qquad\mathsf{b}_{\eta_{k}}^{0}=\mathsf{a}_{\eta_{k}},\qquad\mathsf{c}_{\eta_{k}}^{0}=\mathsf{d}_{\eta_{k}}.
\end{equation}

Moreover, the state associated with the paths are also given by the
same rule as in the periodic case, namely, 
\begin{itemize}
\item The state $\left|P^{(u_{k})}\right\rangle $ associated with any path
$P$ whose the label of its leaf node is $\lambda_{n}^{P}=u_{k}$
is given by, 
\begin{equation}
|P^{(u_{k})}\rangle=\prod_{j=0,j\neq k}^{n}B(u_{j})\Psi_{0}.
\end{equation}
\end{itemize}
The condition for a given path to end with the label $u_{k}$ is that
it contains a filled node (circle or square) in the level $k$ and
no other filled node in the higher levels. Hence, from combinatorial
arguments, it follows that there is only one path in each diagram
ending with the label $u_{0}$ and $2\cdot3^{k-1}$ paths ending with
the label $u_{k}$ for $k\geqslant1$. The partition of the $3^{n}$
terms (the number of paths in these ternary combinatorial tree of
length $n$) into $n+1$ groups (corresponding to paths ending with
the same label) is only possible thanks to the identity $1+2\left(1+3+\cdots+3^{n-1}\right)=3^{n}$.

If we are interested only in the eigenvalues of the transfer matrix
and in the Bethe Ansatz equations, then we shall need only of the
following rules (which are, as a matter of a fact, the same rules
for the periodic case):
\begin{itemize}
\item The eigenvalues are determined by the sum of all paths of the diagrams
$\mathcal{A}$ and $\mathcal{D}$ ending with the label $u_{0}$.
\end{itemize}
Since there is only a path terminating with the label $u_{0}$ in
both the diagrams, we get, promptly, that \begin{equation}
\begin{tikzpicture}
\node(t) at (-3.1,0){$\tau_n\left(u_0|u_1,\ldots,u_n\right)\Psi_n\left(u_1,\ldots,u_n\right)=$};
\node(0) at (0,0) [rectangle, red,thick, draw, inner sep=2.5]{$\mathcal{A}$};
\node(1) at (1,0) [circle,blue   , thick, draw, inner sep=2.5]{};
\node(2) at (2,0) [circle,blue   , thick, draw, inner sep=2.5]{};

\draw[thick] (0)--(1);
\draw[thick,dashed] (1)--(2);

\node(p) at (2.45,0){$+$};
\node(3) at (3,0) [rectangle,red, thick, draw, inner sep=2.5]{$\mathcal{D}$};
\node(4) at (4,0) [circle,blue   , thick, draw, inner sep=2.5]{};
\node(5) at (5,0) [circle,blue   , thick, draw, inner sep=2.5]{};

\draw[thick] (3)--(4);
\draw[thick,dashed] (4)--(5);
   
\end{tikzpicture}
\end{equation} that is, 
\begin{equation}
\tau_{n}\left(u_{0}|u_{1},\ldots,u_{n}\right)=\kappa_{1,1}^{+}\left(u_{0}\right)\mathfrak{a}(u_{0})\prod_{k=1}^{n}\mathsf{a}_{1}\left(u_{0},u_{k}\right)+\kappa_{2,2}^{+}\left(u_{0}\right)\mathfrak{d}(u_{0})\prod_{k=1}^{n}\mathsf{d}_{1}\left(u_{0},u_{k}\right),
\end{equation}
which agrees with (\ref{eigen-B}).

For the boundary Bethe Ansatz equations, we have,
\begin{itemize}
\item The Bethe Ansatz equation fixing the rapidity $u_{k}$ is determined
by the sum of all paths of the $\mathcal{A}$ and $\mathcal{D}$ diagrams
that end with the label $u_{k}$.
\end{itemize}
Therefore, from the values associated with these paths, we get that,
\begin{align}
 & \beta_{n}^{k}\left(u_{0}|u_{1},\ldots,u_{n}\right)\nonumber \\
 & =\kappa_{1,1}^{+}(u_{0})\mathfrak{a}\left(u_{k}\right)\sum_{\eta_{1},\ldots,\eta_{k-1}=1}^{3}\prod_{i=1}^{k-1}\mathsf{b}_{\eta_{i}}^{\eta_{i-1}}\left(\lambda_{i-1}^{P},u_{i}\right)\mathsf{a}_{2}\left(\lambda_{k-1}^{P},u_{k}\right)\prod_{j=k+1}^{n}\mathsf{a}_{1}\left(u_{k},u_{j}\right)\nonumber \\
 & +\kappa_{1,1}^{+}(u_{0})\mathfrak{d}\left(u_{k}\right)\sum_{\eta_{1},\ldots,\eta_{k-1}=1}^{3}\prod_{i=1}^{k-1}\mathsf{b}_{\eta_{i}}^{\eta_{i-1}}\left(\lambda_{i-1}^{P},u_{i}\right)\mathsf{a}_{3}\left(\lambda_{k-1}^{P},u_{k}\right)\prod_{j=k+1}^{n}\mathsf{d}_{1}\left(u_{k},u_{j}\right)\nonumber \\
 & +\kappa_{2,2}^{+}(u_{0})\mathfrak{d}\left(u_{k}\right)\sum_{\eta_{1},\ldots,\eta_{k-1}=1}^{3}\prod_{i=1}^{k-1}\mathsf{c}_{\eta_{i}}^{\eta_{i-1}}\left(\lambda_{i-1}^{P},u_{i}\right)\mathsf{d}_{2}\left(\lambda_{k-1}^{P},u_{k}\right)\prod_{j=k+1}^{n}\mathsf{d}_{1}\left(u_{k},u_{j}\right)\nonumber \\
 & +\kappa_{2,2}^{+}(u_{0})\mathfrak{a}\left(u_{k}\right)\sum_{\eta_{1},\ldots,\eta_{k-1}=1}^{3}\prod_{i=1}^{k-1}\mathsf{c}_{\eta_{i}}^{\eta_{i-1}}\left(\lambda_{i-1}^{P},u_{i}\right)\mathsf{d}_{3}\left(\lambda_{k-1}^{P},u_{k}\right)\prod_{j=k+1}^{n}\mathsf{a}_{1}\left(u_{k},u_{j}\right)=0, &  & 1\leqslant k\leqslant n,
\end{align}
or, grouping the terms containing $\mathfrak{a}\left(u_{k}\right)$
and $\mathfrak{d}\left(u_{k}\right)$, 
\begin{align}
 & \beta_{n}^{k}\left(u_{0}|u_{1},\ldots,u_{n}\right)=\mathfrak{a}\left(u_{k}\right)\prod_{j=k+1}^{n}\mathsf{a}_{1}\left(u_{k},u_{j}\right)\times\nonumber \\
 & \sum_{\eta_{1},\ldots,\eta_{k-1}=1}^{3}\left[\prod_{i=1}^{k-1}\mathsf{b}_{\eta_{i}}^{\eta_{i-1}}\left(\lambda_{i-1}^{P},u_{i}\right)\mathsf{a}_{2}\left(\lambda_{k-1}^{P},u_{k}\right)\kappa_{1,1}^{+}(u_{0})+\prod_{i=1}^{k-1}\mathsf{c}_{\eta_{i}}^{\eta_{i-1}}\left(\lambda_{i-1}^{P},u_{i}\right)\mathsf{d}_{3}\left(\lambda_{k-1}^{P},u_{k}\right)\kappa_{2,2}^{+}(u_{0})\right]\nonumber \\
 & +\mathfrak{d}\left(u_{k}\right)\prod_{j=k+1}^{n}\mathsf{d}_{1}\left(u_{k},u_{j}\right)\times\nonumber \\
 & \sum_{\eta_{1},\ldots,\eta_{k-1}=1}^{3}\left[\prod_{i=1}^{k-1}\mathsf{c}_{\eta_{i}}^{\eta_{i-1}}\left(\lambda_{i-1}^{P},u_{i}\right)\mathsf{d}_{2}\left(\lambda_{k-1}^{P},u_{k}\right)\kappa_{2,2}^{+}(u_{0})+\prod_{i=1}^{k-1}\mathsf{b}_{\eta_{i}}^{\eta_{i-1}}\left(\lambda_{i-1}^{P},u_{i}\right)\mathsf{a}_{3}\left(\lambda_{k-1}^{P},u_{k}\right)\kappa_{1,1}^{+}(u_{0})\right],\nonumber \\
 & 1\leq k\leq n.\label{BetheBold}
\end{align}

The simplest Bethe Ansatz equation is that one fixing the rapidity
$u_{1}$. Differently from the periodic case, we have here two paths
ending with the label $u_{1}$ in each diagram \textendash{} namely,
the paths $P\left(2,1,\ldots,1\right)$ and $P\left(3,1,\ldots,1\right)$.
Notice moreover that the leaf node of the path $P_{\mathcal{A}}\left(2,1,\ldots,1\right)$
is a circle, while the leaf node of $P_{\mathcal{A}}\left(3,1,\ldots,1\right)$
is a square (similarly, the leaf nodes in the paths $P_{\mathcal{D}}\left(2,1,\ldots,1\right)$
and $P_{\mathcal{D}}\left(3,1,\ldots,1\right)$ are respectively a
square and a circle). From this we get that \begin{equation}
\begin{tikzpicture}
\node(t) at (0,0)[left]{$\beta_n^1\left(u_0|u_1,\ldots,u_n\right)\Psi_n\left(u_{j}^{\times}\right)=\hspace{0.2cm}$};
\node(00) at (0,0) [rectangle,red, thick, draw, inner sep=2.5]{$\mathcal{A}$};
\node(11) at (1, 0.5) [circle,blue   ,fill=blue, thick, draw, inner sep=2.5]{};
\node(21) at (2, 0.5) [circle,blue   , thick, draw, inner sep=2.5]{};
\node(31) at (3, 0.5) [circle,blue   , thick, draw, inner sep=2.5]{};
\node(12) at (1,-0.5) [rectangle,blue   ,fill=blue, thick, draw, inner sep=3]{};
\node(22) at (2,-0.5) [rectangle,blue   , thick, draw, inner sep=3]{};
\node(32) at (3,-0.5) [rectangle,blue   , thick, draw, inner sep=3]{};

\draw[thick,cap=round] (00) -- (11);
\draw[thick,cap=round] (11) -- (21);
\draw[thick,cap=round] (21) -- (31);
\draw[thick,cap=round] (00) -- (12);
\draw[thick,cap=round] (12) -- (22);
\draw[thick,cap=round] (22) -- (32);

\node(P)  at (3.45,0){$+$};

\node(41) at (4, 0  ) [rectangle,red, thick, draw, inner sep=2.5]{$\mathcal{D}$};
\node(51) at (5, 0.5) [rectangle,blue,fill=blue, thick, draw, inner sep=3]{};
\node(61) at (6, 0.5) [rectangle,blue, thick, draw, inner sep=3]{};
\node(71) at (7, 0.5) [rectangle,blue, thick, draw, inner sep=3]{};
\node(52) at (5,-0.5) [circle,blue   ,fill=blue, thick, draw, inner sep=2.5]{};
\node(62) at (6,-0.5) [circle,blue   , thick, draw, inner sep=2.5]{};
\node(72) at (7,-0.5) [circle,blue   , thick, draw, inner sep=2.5]{};

\draw[thick,cap=round] (41) -- (51);
\draw[thick,cap=round] (51) -- (61);
\draw[thick,cap=round] (61) -- (71);
\draw[thick,cap=round] (41) -- (52);
\draw[thick,cap=round] (52) -- (62);
\draw[thick,cap=round] (62) -- (72);

\end{tikzpicture}
\end{equation}that is, 
\begin{multline}
\beta_{n}^{1}\left(u_{0}|u_{1},\ldots,u_{n}\right)=\mathsf{a}_{2}\left(u_{0},u_{1}\right)\prod_{i=2}^{n}\mathsf{a}_{1}\left(u_{1},u_{i}\right)\kappa_{1,1}^{+}\left(u_{0}\right)\mathfrak{a}\left(u_{1}\right)+\mathsf{a}_{3}\left(u_{0},u_{1}\right)\prod_{i=2}^{n}\mathsf{d}_{1}\left(u_{1},u_{i}\right)\kappa_{1,1}^{+}\left(u_{0}\right)\mathfrak{d}\left(u_{1}\right)\\
+\mathsf{d}_{2}\left(u_{0},u_{1}\right)\prod_{i=2}^{n}\mathsf{d}_{1}\left(u_{1},u_{i}\right)\kappa_{2,2}^{+}\left(u_{0}\right)\mathfrak{d}\left(u_{1}\right)+\mathsf{d}_{3}\left(u_{0},u_{1}\right)\prod_{i=2}^{n}\mathsf{a}_{1}\left(u_{1},u_{i}\right)\kappa_{2,2}^{+}\left(u_{0}\right)\mathfrak{a}\left(u_{1}\right),
\end{multline}
 or, grouping the common terms, 
\begin{multline}
\beta_{n}^{1}\left(u_{0}|u_{1},\ldots,u_{n}\right)=\left[\mathsf{a}_{2}\left(u_{0},u_{1}\right)\kappa_{1,1}^{+}\left(u_{0}\right)+\mathsf{d}_{3}\left(u_{0},u_{1}\right)\kappa_{2,2}^{+}\left(u_{0}\right)\right]\mathfrak{a}\left(u_{1}\right)\prod_{i=2}^{n}\mathsf{a}_{1}\left(u_{1},u_{i}\right)\\
+\left[\mathsf{d}_{2}\left(u_{0},u_{1}\right)\kappa_{2,2}^{+}\left(u_{0}\right)+\mathsf{a}_{3}\left(u_{0},u_{1}\right)\kappa_{1,1}^{+}\left(u_{0}\right)\right]\mathfrak{d}\left(u_{1}\right)\prod_{i=2}^{n}\mathsf{d}_{1}\left(u_{1},u_{i}\right)=0.
\end{multline}
Finally, thanks to the symmetry of the excited states regarding the
permutation of the rapidities, the other Bethe Ansatz equations given
at (\ref{BetheBold}) can be simplified, as we did in the periodic
case. In fact, (\ref{BB2}) implies that both the $\mathcal{A}$ as
well as the $\mathcal{D}$ diagrams must be symmetric with respect
to the permutation of their labels. From this we conclude that the
Bethe Ansatz equation fixing the rapidity $u_{k}$ can actually be
written in the same form as that fixing $u_{1}$, except that the
rapidities $u_{k}$ and $u_{1}$ must be permuted. That is, likewise
in the periodic case, we have the rule, 
\begin{itemize}
\item The Bethe Ansatz equation fixing the rapidity $u_{k}$ $(2\leqslant k\leqslant n)$
can be obtained by the same paths that determine the Bethe Ansatz
equation fixing $u_{1}$, provided that the labels $u_{1}$ and $u_{k}$
are permuted.
\end{itemize}
This leads us to the final form of the Bethe Ansatz equations for
the six-vertex model with boundaries: 
\begin{multline}
\beta_{n}^{k}\left(u_{0}|u_{1},\ldots,u_{n}\right)=\left[\mathsf{a}_{2}\left(u_{0},u_{k}\right)\kappa_{1,1}^{+}\left(u_{0}\right)+\mathsf{d}_{3}\left(u_{0},u_{k}\right)\kappa_{2,2}^{+}\left(u_{0}\right)\right]\mathfrak{a}\left(u_{k}\right)\prod_{i=1,i\neq k}^{n}\mathsf{a}_{1}\left(u_{k},u_{i}\right)\\
+\left[\mathsf{d}_{2}\left(u_{0},u_{k}\right)\kappa_{2,2}^{+}\left(u_{0}\right)+\mathsf{a}_{3}\left(u_{0},u_{k}\right)\kappa_{1,1}^{+}\left(u_{0}\right)\right]\mathfrak{d}\left(u_{k}\right)\prod_{i=1,i\neq k}^{n}\mathsf{d}_{1}\left(u_{k},u_{i}\right)=0,\label{BAEBfinal}
\end{multline}
 which is in agreement with (\ref{BAE-B}).

Finally, the comparison between (\ref{BetheBold}) and (\ref{BAEBfinal})
provides us as well with other more intricate identities:
\begin{align}
\sum_{\eta_{1},\ldots,\eta_{k-1}=1}^{3}\prod_{i=1}^{k-1}\mathsf{b}_{\eta_{i}}^{\eta_{i-1}}\left(\lambda_{i-1}^{P},u_{i}\right)\mathsf{a}_{2}\left(\lambda_{k-1}^{P},u_{k}\right) & =\mathsf{a}_{2}\left(u_{0},u_{k}\right)\prod_{i=1}^{k-1}\mathsf{a}_{1}\left(u_{k},u_{i}\right),\\
\sum_{\eta_{1},\ldots,\eta_{k-1}=1}^{3}\prod_{i=1}^{k-1}\mathsf{b}_{\eta_{i}}^{\eta_{i-1}}\left(\lambda_{i-1}^{P},u_{i}\right)\mathsf{a}_{3}\left(\lambda_{k-1}^{P},u_{k}\right) & =\mathsf{a}_{3}\left(u_{0},u_{k}\right)\prod_{i=1}^{k-1}\mathsf{d}_{1}\left(u_{k},u_{i}\right),\\
\sum_{\eta_{1},\ldots,\eta_{k-1}=1}^{3}\prod_{i=1}^{k-1}\mathsf{c}_{\eta_{i}}^{\eta_{i-1}}\left(\lambda_{i-1}^{P},u_{i}\right)\mathsf{d}_{2}\left(\lambda_{k-1}^{P},u_{k}\right) & =\mathsf{d}_{2}\left(u_{0},u_{k}\right)\prod_{i=1}^{k-1}\mathsf{d}_{1}\left(u_{k},u_{i}\right),\\
\sum_{\eta_{1},\ldots,\eta_{k-1}=1}^{3}\prod_{i=1}^{k-1}\mathsf{c}_{\eta_{i}}^{\eta_{i-1}}\left(\lambda_{i-1}^{P},u_{i}\right)\mathsf{d}_{3}\left(\lambda_{k-1}^{P},u_{k}\right) & =\mathsf{d}_{3}\left(u_{0},u_{k}\right)\prod_{i=1}^{k-1}\mathsf{a}_{1}\left(u_{k},u_{i}\right).
\end{align}

\section{Conclusion\label{Section-Conclusion}}

In this paper we presented a combinatorial approach for the (periodic
and non-periodic) algebraic Bethe Ansatz in which the commutation
relations are represented by combinatorial diagrams, so that the action
the diagonal operators on the $n$th excited states becomes described
by labeled combinatorial trees. From the analysis of these combinatorial
diagrams, every term resulting from the algebraic Bethe Ansatz can
be easily recovered following simple rules. In special, the eigenvalues
are obtained in a straightforward way, namely, they are provided by
those paths containing only hollow nodes. The Bethe Ansatz equations
also can be found directly from this analysis: we can first obtain
the Bethe Ansatz Ansatz equation fixing the rapidity $u_{1}$, which
is provided by the paths containing a filled node in the first level
and no other filled node beyond that, and then, thanks to the symmetry
of the diagrams regarding the permutation of their labels, we get
that the Bethe Ansatz equations fixing the rapidity $u_{k}$ will
have the same form of the previous one, except that the rapidities
$u_{1}$ and $u_{k}$ should be permuted. This symmetry also provides
several mathematical intricate identities. 

As a first generalization of this combinatorial approach could be
the analysis of the fifteen and nineteen vertex models. These cases
are interesting because the creator operators do not commute among
themselves anymore, which leads to a more elaborate form for the excited
states. We expect nevertheless that these excited states can also
be determined through the analysis of some combinatorial trees, for
instance, imposing the symmetry requirements regarding the permutation
of their labels. In this way, the present analysis could provide the
eigenvalues and the eigenstates of the transfer matrix, and also the
corresponding Bethe Ansatz equations, in a more direct way. Besides,
perhaps this combinatorial method can be useful as well in the analysis
of scalar products and correlation functions of integrable models.

Finally, since we also presented a comprehensive introduction to the
algebraic Bethe Ansatz of the six-vertex model \textendash{} for both
periodic and non-periodic boundary conditions \textendash , we believe
that this work can also be useful for teaching the algebraic Bethe
Ansatz for introductory audiences.
\begin{acknowledgments}
This work was supported in part by São Paulo Research Foundation (FAPESP,
grant \#2011/18729-1), Brazilian Research Council (CNPq, grant \#310625/2013-0),
and Coordination for the Improvement of Higher Education Personnel
(CAPES).
\end{acknowledgments}

\bibliographystyle{aipnum4-1}
\bibliography{Trees}

\end{document}